\documentclass{article}
\usepackage{latexsym}
\usepackage{amsmath,amssymb,amsthm}

\newtheorem{theorem}{Theorem}[section]
\newtheorem{corollary}[theorem]{Corollary}

\newtheorem{lemma}[theorem]{Lemma}
\newtheorem{remark}[theorem]{Remark}
 \makeatletter
\@addtoreset{equation}{section}

\makeatother
 
\newcommand{\app}[4]{F_{\!#1}\!
  \left(\left.{#2 \atop #3}\right| #4 \right) }

\newcommand{\hpg}[5]{{}_{#1}\mbox{\rm F}_{\!#2}\!
  \left(\left.{#3 \atop #4}\right| #5 \right) }
\newcommand{\hpgg}[4]{{}_{#1}\mbox{\rm F}_{\!#2}\!
  \left({#3 \atop #4}\right) }
\newcommand{\hpgo}[2]{{}_{#1}\mbox{\rm F}_{\!#2}}

\newcommand{\equal}{&\!\!\!=\!\!\!&}

\newcommand{\cc}{\lambda}

\newcommand{\CC}{{\Bbb C}}
\newcommand{\QQ}{{\Bbb Q}}
\newcommand{\PP}{{\Bbb P}}
\newcommand{\ZZ}{{\Bbb Z}}

\title
{Transformations and invariants for dihedral Gauss hypergeometric functions}

\author{Raimundas Vidunas\\
\em Kobe University}

\begin{document}

\maketitle

\begin{abstract} 
Hypergeometric equations with a dihedral monodromy group can be solved in terms
of elementary functions. This paper gives explicit general expressions for quadratic 
monodromy invariants for these hypergeometric equations, using a generalization
of Clausen's formula and terminating double hypergeometric sums.
Besides, pull-back transformations for the dihedral hypergeometric equations
are presented, including Klein's pullback transformations for 
the equations with a finite (dihedral) monodromy group. 
\end{abstract}

\section{Introduction}
\label{logarithms}

By a {\em dihedral hypergeometric equation} we mean a hypergeometric equation
\begin{equation} \label{hpgde}
z\,(1-z)\,\frac{d^2y(z)}{dz^2}+
\big(C-(A+B+1)\,z\big)\,\frac{dy(z)}{dz}-A\,B\,y(z)=0
\end{equation}
with a dihedral monodromy group. The hypergeometric solutions 
of these equations are called {\em dihedral Gauss hypergeometric functions}, 
or {\em dihedral $\hpgo21$ functions}. These functions form a class
of hypergeometric functions that can be represented  in terms of elementary functions:
power and rational functions, and logarithms or inverse trigonometric functions.
They are {\em Liouvillian functions}  \cite{singulm1} since the dihedral monodromy group is solvable.  
General elementary expressions for the dihedral functions are presented in \cite{tdihedral} 
and recalled in Subsection \ref{sec:explicit} here.

The main purpose of this paper is to give elementry expressions for monodromy invariants
for the dihedral hypergeometric equations. A general expression for quadratic monodromy
invariants is obtained from a generalization \cite{gclausen} of Claussen's formula
\begin{equation} \label{eq:clausen}
\hpg21{A,\;B}{A+B+\frac12}{\,x}^2=\hpg32{2A,\,2B,\,A+B}{2A+2B,A+B+\frac12}{\,x}.
\end{equation}
The elementary expressions involve terminating double hypergeometric sums;
they are presented in Section \ref{sec:symsq}. 
In the cases of finite dihedral monodromy group of order $2n$, 
the hypergeometric solutions are algebraic functions and there are additionally
two independent invariants of degree $n$. A general closed form of these invariants is
desirable \cite{yoshidih}, \cite{ochiay}. This paper gives a simple
recipe to compute these invariants in Section \ref{adihedral},
and differential expressions for them in Subsection \ref{diffinv}.

A related topic is pull-back transformations of dihedral hypergeometric equations,
and related algebraic transformations of dihedral hypergeometric functions.
A {\em pull-back transformation} of a differential equation for $y(z)$ in $d/dz$ 
has the form
\begin{equation} \label{algtransf}
z\longmapsto\varphi(x), \qquad y(z)\longmapsto
Y(x)=\theta(x)\,y(\varphi(x)),
\end{equation}
where $\varphi(x)$ is a rational function, and  $\theta(x)$ is
a product of powers of rational functions.
Geometrically, the transformation {\em pull-backs} the starting
differential equation on the projective line $\PP^1_z$ to a
differential equation on the projective line $\PP^1_x$, with
respect to the covering $\varphi:\PP^1_x\to\PP^1_z$
determined by the rational function $\varphi(x)$.
Fuchsian equations are pull-backed to Fuchsian equations again,
usually with more singular points. 
The factor $\theta(x)$ shifts the local exponents of the pull-backed equation,
but it does not change the local exponent differences.

Of particular interest are pull-back transformations of hypergeometric equation (\ref{hpgde})
into itself, perhaps with other values of the parameters $A,B,C$. The pull-back transformations 
induce transformation formulas for the Gauss hypergeometric function of the form
\begin{equation} \label{hpgtransf}
\hpg{2}{1}{\!\widetilde{A},\,\widetilde{B}\,}{\widetilde{C}}{\,x}
=\theta(x)\;\hpg{2}{1}{\!A,\,B}{C}{\varphi(x)}.
\end{equation}
These transformations are classified in \cite{algtgauss}. It was widely assumed that the
classical quadratic  and degree 3, 4, 6 transformations (due to Gauss, Kummer, Goursat) 
exhaust all hypergeometric transformations of this form. 
Particularly, \cite[Section 2.1.5]{bateman} states the following: 
\begin{quote} 
Transformations of [degrees other than $2, 3, 4, 6$] can exist
only if $a,b,c$ are certain rational numbers; in those cases the solutions
of the hypergeometric equation are algebraic functions."\\
\end{quote}
Transformations of general dihedral hypergeometric functions (of any degree) provide
a simple counter-example to this assumption, as noticed in \cite{algtgauss}, \cite{andkitaev}. 
The transformations with a continuous parameter are presented here
exhaustively in Section \ref{sec:trihedr}. Furthermore, the well-known theorem of Klein 
implies that any hypergeometric equation with a finite dihedral monodromy group is 
a pullback of the simplest hypergeometric equation with the same monodromy group.
These transformations are considered in Sections \ref{adihedral}. 

\section{Preliminary facts}

Because of frequent use, we recall Euler's and Pfaff's
fractional-linear transformations \cite[Theorem 2.2.5]{specfaar}:
\begin{eqnarray} \label{flinear1}
\hpg{2}{1}{A,\,B\,}{C}{\,z} \equal
(1-z)^{C-A-B}\;\hpg{2}{1}{C-A,\,C-B}{C}{\,z} \\
\label{flinear3} \equal (1-z)^{-A}\;\hpg{2}{1}{A,\,C-B\,}{C}{\frac{z}{z-1}}.
\end{eqnarray}
As usual, $(a)_n$ denotes the {\em Pochhammer symbol}  (also called {\em raising factorial}),
which is the product $a(a+1)\ldots(a+n-1)$. The {\em double factorial} is used to adjust the look
of some of our formulas:
\begin{equation} \label{eq:dfact}
(2k-1)!!=2^k\left({\textstyle \frac12}\right)_k=\frac{(2k)!}{2^{k}\,k!}, \qquad 
(2k+1)!!=2^{k+1}\left({\textstyle \frac12}\right)_{k+1}=\frac{(2k+1)!}{2^{k}\,k!},
\end{equation}
and $(2k)!!=2^kk!$.

\subsection{Fuchsian differential equations}

Hypergeometric equation (\ref{hpgde}) is 
a canonical Fuchsian equation on $\PP^1$ 
with three singular points. The singularities are $z=0,1,\infty$, 
and the local exponents are:
\[
\mbox{$0$, $1-C$ at $z=0$;} \qquad \mbox{$0$, $C-A-B$ at $z=1$;}\qquad
\mbox{and $A$, $B$ at $z=\infty$.}
\]
The local exponent differences at the singular points are equal (up to a sign) to $1-C$,
\mbox{$C-A-B$} and $A-B$, respectively. We denote a hypergeometric equation 
with the local exponent differences $d_1,d_2,d_3$ by $E(d_1,d_2,d_3)$, 
and consider the order of the three arguments unimportant.

For the dihedral hypergeometric equations, two of the local exponent differences
are half-integers. The dihedral hypergeometric functions are  {\em contiguous}
to the hypergeometric equation $E(1/2,1/2,a)$ with $a\in\CC$.

A canonical Fuchsian equation with 4 singular points is Heun's equation
\begin{align}  \label{eq:heun}
\frac{d^2y(x)}{dx^2}+\biggl(\frac{\gamma}{x}+\frac{\delta}{x-1}
+\frac{\alpha+\beta-\gamma-\delta+1}{x-t}\biggr)\frac{dy(x)}{dx}
+\frac{\alpha\beta x-q}{x(x-1)(x-t)} y(x) = 0.
\end{align}
The local exponents are:
\[
\mbox{$0$, $1-\gamma$ at $x=0$;} \quad 
\mbox{$0$, $1-\delta$ at $x=1$;} \quad 
\mbox{$0$, $\gamma+\delta-\alpha-\beta$ at $x=t$;}\quad
\mbox{$\alpha$, $\beta$ at $x=\infty$.}
\]
The coefficients of power series solutions of Heun's equation satisfy a rather complicated 
second order recurrence \cite{MaierHeun}.
The quadratic substitution $z\mapsto 4x(1-x)$ transforms hypergeometric equation 
(\ref{hpgde}) to Heun's equation with
\begin{equation} \textstyle
(\alpha,\,\beta,\,\gamma,\,\delta,\,t,\,q)=\left(2A,\,2B,\,C,\,C,\,\frac12,\,2AB\right). 
\end{equation}
Other example of a quadratic transformation between hypergeometric and Heun functions
is given by formulas (\ref{heunhpg2a})--(\ref{heunhpg2b}) below.

\subsection{Double hypergeometric series}

Quadratic invariants of dihedral hypergeometric equations are expressed in terms
of the following double hypergeometric series:
\begin{eqnarray} \label{eq:gap211}
F^{2:1;1}_{1:1;1}\!\left(\left. {a; b; u_1, u_2\atop c; v_1, v_2} \right| x,y\right) \!\equal
\sum_{p=0}^{\infty} \sum_{q=0}^{\infty}
\frac{(a)_{p+q}\,(b)_{p+q}\,(u_1)_p\,(u_2)_q}{(c)_{p+q}\,(v_1)_p\,(v_2)_q\;p!\,q!}\,x^p\,y^q,\\
\label{eq:gap122}
F^{1:2;2}_{2:0;0}\!\left(\left. {a; u_1, u_2; v_1, v_2 \atop b; \; c} \right| x,y\right) \!\equal
\sum_{p=0}^{\infty} \sum_{q=0}^{\infty}
\frac{(a)_{p+q}\,(u_1)_p\,(u_2)_q\,(v_1)_p\,(v_2)_q}{(b)_{p+q}\;(c)_{p+q}\;p!\,q!}\,x^p\,y^q.\qquad
\end{eqnarray}
If $a=c$, these two series reduce to Appell's bivariate hypergeometric series 
$$\app2{b; u_1,u_2}{v_1,v_2}{x,y},\qquad 
\app3{u_1,u_2;v_1,v_2}{b}{x,y},$$
respectively. We use only terminating versions of these double hypergeometric sums.

Appell's $F_2(x,y)$ and $F_3(x,y)$ functions are closely related. 
They satisfy the same system of partial differential equations up to a simple transformation,
and terminating $F_2$ sums become terminating $F_3$ sums when 
summation is reversed in both directions. In particular, for $(a)_{k+\ell}\neq 0$ we have
\begin{eqnarray} \label{eq:r2f3rel}
\app2{a; -k,-\ell}{-2k,-2\ell}{x,y}=\frac{k!\,\ell!\,(a)_{k+\ell}}{(2k)!\,(2\ell)!}\,x^k\,y^{\ell}\,
\app3{k+1,\ell+1; -k,-\ell}{1-a-k-\ell}{\frac1x,\frac1y}.
\end{eqnarray}
The relation between $F^{2:1;1}_{1:1;1}(x,y)$ or $F^{1:2;2}_{2:0;0}(x,y)$ hypergeometric functions
is similar \cite{gclausen}. In particular, we will refer to the following formula for reversing summation order in both directions in a terminating $F^{2:1;1}_{1:1;1}$ sum:
\begin{eqnarray} \label{eq:reverse32}
\hspace{-5pt} F^{2:1;1}_{1:1;1}\!\left(\left. {\!a; -a-2k-2\ell; -k, -\ell\atop \frac12-k-\ell; -2k, -2\ell} 
\right| x,1-x\right) = \hspace{158pt} \nonumber \\ 
 \frac{(a)_{2k+2\ell+1}\,k!\,\ell!\,x^k\,(1-x)^\ell}
{(a+k+\ell)\left(\frac12\right)_{k+\ell}(2k)!\,(2\ell)!}\, 
F^{1:2;2}_{2:0;0}\!\left(\left. {\frac12;\, k+1,\ell+1; -k, -\ell\atop a+k+\ell+1;1-a-k-\ell} 
\right|  \frac{1}{x},\frac1{1-x}\right). 
\end{eqnarray}

Note that hypergeometric series like $\displaystyle\app2{a; -k,-\ell}{-2k,-2\ell}{x,y}$ is not conventionally
defined for a non-negative integers $k,\ell$, because of the zero or negative lower parameter.
As described in \cite[Remark 3.2]{tdihedral}, we choose a terminating interpretation
of these hypergeometric series.

\subsection{Explicit expressions for dihedral functions}
\label{sec:explicit}

The following formulas are elementary expressions for general dihedral hypergeometric functions
proved in \cite{tdihedral}. The numbers $k,\ell$ are non-negative integers, and $a\in\CC\setminus\ZZ$.
The $F_2$ and $F_3$ series are 
finite sums of $(k+1)(\ell+1)$ terms. 
\begin{eqnarray} 
\label{eq:diha} \hspace{-9pt}
\hpg{2}{1}{\frac{a}{2},\,\frac{a+1}{2}+\ell}{a+k+\ell+1}{1-z}  \!\equal\!
z^{k/2} \left(\frac{1+\sqrt{z}}{2} \right)^{-a-k-\ell}
\times\nonumber\\
&& \! \app3{k+1,\ell+1; -k,-\ell}{a+k+\ell+1}{\frac{\sqrt{z}-1}{2\sqrt{z}},\frac{1-\sqrt{z}}2},\\
\label{eq:dih12} \hspace{-9pt}
\frac{\left(\frac{a+1}2\right)_\ell}{\left(\frac12\right)_\ell}\,
\hpg{2}{1}{\frac{a}{2},\frac{a+1}{2}+\ell}{\frac{1}{2}-k}{\,z} \!\equal\!
\frac{(1+\sqrt{z})^{-a}}2
\app2{a; -k,-\ell}{-2k,-2\ell}{\frac{2\sqrt{z}}{1+\sqrt{z}},\frac2{1+\sqrt{z}}}\nonumber\\
&& \hspace{0pt}+\frac{(1-\sqrt{z})^{-a}}2
\app2{a; -k,-\ell}{-2k,-2\ell}{\frac{2\sqrt{z}}{\sqrt{z}-1},\frac2{1-\sqrt{z}}},\\ 
\label{eq:dih32} \hspace{-11pt}
\frac{\left(\frac{a+1}2\right)_k\left(\frac{a}2\right)_{k+\ell+1}}
{\left(\frac12\right)_k\left(\frac12\right)_{k+1}\left(\frac12\right)_\ell}\,(-1)^k z^{k+\frac12}\,
\hpg{2}{1}{\frac{a+1}{2}+k,\,\frac{a}{2}+k+\ell+1\,}{\frac{3}{2}+k}{\,z}\hspace{-152pt}\nonumber\\
\!\equal\! \frac{(1-\sqrt{z})^{-a}}2
\app2{a; -k,-\ell}{-2k,-2\ell}{\frac{2\sqrt{z}}{\sqrt{z}-1},\frac2{1-\sqrt{z}}}\nonumber\\
&&\hspace{0pt}-\frac{(1+\sqrt{z})^{-a}}2
\app2{a; -k,-\ell}{-2k,-2\ell}{\frac{2\sqrt{z}}{1+\sqrt{z}},\frac2{1+\sqrt{z}}}. \qquad
\end{eqnarray}
The $F_2$ and $F_3$ sums become univariate $\hpgo21$ sums when $k=0$ or $\ell=0$,
and they become the constant 1 when $k=\ell=0$. The simplest expressions for the $k=\ell=0$
case are well-known \cite[15.1]{abrostegun}, \cite[2.8]{bateman}. 
In particular,
\begin{eqnarray}  \label{dihedr2}
\hpg{2}{1}{\frac{a}{2},\,\frac{a+1}{2}\,}{\frac{1}{2}}{\,z} \equal
\frac{(1-\sqrt{z})^{-a}+(1+\sqrt{z})^{-a}}{2}.
\end{eqnarray}
If $k=\ell$, then the terminating $F_2$ of $F_3$ sums can be likewise replaced
by terminating $\hpgo21$ sums, since this case reduces to $k\ell=0$ because
of a quadratic transformation. For example, formula \cite[(3.1.9)]{specfaar} gives 
\begin{equation} \label{eq:quadrkk}
\hpg21{a,a+k+\frac12}{\frac12-k}{z} =
(1+z)^{-a}\,\hpg21{\frac{a}2,\,\frac{a+1}2}{\frac12-k}{\frac{4z}{(1+z)^2}}.
\end{equation}
The trigonometric substitution $z\mapsto -\tan^2 x$ gives attractive trigonometric formulas
for the dihedral functions \cite[Section 6]{tdihedral}. 
General formulas for the logarithmic and degenerate cases with $a\in\ZZ$
are given in \cite{tdihedral} as well. 

The proof of general expressions (\ref{eq:diha})--(\ref{eq:dih32}) 
is based on the fact \cite{AppelGhpg} that the univariate functions
\begin{equation}\label{eq:f2gauss}
\app2{a;\;b_1,b_2}{2b_1,2b_2}{x,2-x}  \quad\mbox{and}\quad 
(x-2)^{-a}\,\hpg21{\frac{a}2,\frac{a+1}2-b_2}{b_1+\frac12}{\frac{x^2}{(2-x)^2}}
\end{equation}
satisfy the same second order Fuchsian equation for general $a,b_1,b_2$.
Next to (\ref{eq:r2f3rel}), we have the following symmetries of
relevant $F_2$ and $F_3$ functions by \cite[Lemma 3.2]{tdihedral}:
\begin{eqnarray}  \label{eq:f2f2rel}
&& \hspace{-20pt} \left(1+\sqrt{z}\right)^{k+\ell}
\app2{a;\,-k,-\ell}{-2k,-2\ell}{\frac{2\sqrt{z}}{1+\sqrt{z}},\frac2{1+\sqrt{z}}}=\nonumber\\
 &&
 \frac{(-1)^\ell\left(\frac{a+1}2\right)_\ell}{\left(\frac{a+1}2+k\right)_\ell}
\left(1-\sqrt{z}\right)^{k+\ell}
\app2{-a-2k-2\ell;\,-k,-\ell}{-2k,-2\ell}{\frac{2\sqrt{z}}{\sqrt{z}-1},\frac2{1-\sqrt{z}}}.\qquad \\
\label{eq:f3f3rel} && \hspace{-20pt}
\app3{k+1,\ell+1; -k,-\ell}{a+k+\ell+1}{\frac{\sqrt{z}-1}{2\sqrt{z}},\frac{1-\sqrt{z}}2}=\nonumber\\
&&\
\frac{(a)_{k+\ell}\left(\frac{a+1}2+k\right)_\ell}{(1+a+k+\ell)_{k+\ell}\left(\frac{a+1}2\right)_\ell}
\,\app3{k+1,\ell+1; -k,-\ell}{1-a-k-\ell}{\frac{\sqrt{z}+1}{2\sqrt{z}},\frac{1+\sqrt{z}}2}.
\end{eqnarray}
The local monodromy differences of the hypergeometric equations
for the dihedral functions here are 
$k+1/2$ $\ell+1/2$ and $a+k+\ell.$ 

\subsection{Quadratic pull-back transformations}

Classical quadratic transformations \cite[Section 3.1]{specfaar} of Gauss hypergeometric functions
illustrate the pull-back transformation
\begin{equation} \label{eq:quadrtr} \textstyle
E\left( \frac12,\cc,\mu \right) \stackrel2\longleftarrow E\left( \cc,\cc,2\mu \right),
\end{equation}
of hypergeometric equations. The arrow notation is the same as in \cite{algtgauss}:
it follows the direction of the pull-back covering $\varphi:\PP^1_x\to\PP^1_z$
and indicates the covering degree 2. 
The local exponent differences double at the two branching points of the quadratic covering, 
and the $\PP^1_z$ point with the local exponent difference $1/2$ is transformed to
a non-singular point. In particular, the quadratic transformation (\ref{eq:quadrkk})
between the $k=\ell$ and $k\ell=0$ cases is
\begin{eqnarray} \label{eq:quadih}  \textstyle
E\left(\frac12,\,k+\frac12,\,\cc\right) \stackrel2\longleftarrow E\left(k+\frac12,\,k+\frac12,\,2\cc \right),
\qquad \mbox{with $\lambda=a+k$.}
\end{eqnarray}

A general hypergeometric equation is pull-backed by the same quadratic transformation 
to a Fuchsian equation with (at least) 4 singular points. 
The local exponent differences transform as
\begin{equation} \textstyle
\left( \cc,\mu,\xi \right) \stackrel2\longleftarrow \left( \cc,\cc,2\mu,2\xi \right),
\end{equation}
and generally Heun's equation is obtained after the transformation.
If the quadratic transformation is applied to a dihedral hypergeometric equation
and the covering ramifies over the two singular points with
half-integer local exponent differences, the local exponent differences transform as 
\begin{equation} \label{eq:qled} \textstyle
\left( k+\frac12,\ell+\frac12,\cc \right) \stackrel2\longleftarrow \left( 2k+1,2\ell+1,\cc,\cc \right).
\end{equation}
After an appropriate choice of the power factor $\theta(x)$ in (\ref{algtransf}), 
the pull-backed equation has trivial monodromy around the two points 
with the integer local exponent differences $2k+1$, $2\ell+1$, 
the global monodromy group is cyclic, 
and the monodromy representation is reducible. 
Properly normalized, the pull-backed Heun equation must have polynomial solutions.
If we apply pull-back transformation (\ref{algtransf}) with 
$\varphi(x)=4x/(1+x)^2$, $\theta(x)=(1+x)^{-a}$ like in (\ref{eq:quadrkk})
to the hypergeometric equation for 
\begin{equation} \label{heunhpg2a}
\displaystyle\hpg{2}{1}{\frac{a}{2},\,\frac{a+1}{2}+\ell}{a+k+\ell+1}{z},
\end{equation}
 the transformed Heun equation has
\begin{equation} \label{heunhpg2b}
(\alpha,\,\beta,\,\gamma,\,\delta,\,t,\,q)=\big(a,\,-k-\ell,\,1+a+k+\ell,\,-2k,\,-1,\,(k-\ell)a \big).
\end{equation}
As demonstrated by formulas (\ref{eq:diha})--(\ref{eq:dih32}), solutions of this
Heun equation can be expressed in terms of terminating Appell's $F_2$ or $F_3$ series,
with the square root of $z$ gone after the substitution $z=1-\varphi(x)$. 
If $\ell=0$ (or $k=0$), the transformation is to a hypergeometric equation
but with a cyclic monodromy group:
\begin{eqnarray} \label{eq:quatriv}  \textstyle
E\left( \frac12,k+\frac12,\cc \right) \stackrel2\longleftarrow E\left( 2k+1,\cc,\cc \right).
\end{eqnarray}
If both $k=0$, $\ell=0$, the transformed equation has just two singularities.

\section{Symmetric square solutions}
\label{sec:symsq}

Monodromy invariants of a hypergeometric equation with a dihedral of finite monodromy group
are elementary solutions of proper symmetric tensor powers of the hypergeometric equation 
\cite{singulm1}. For a general dihedral equation, the tensor square has to be considered. 

Let $y_1,y_2$ denote two independent solutions of hypergeometric equation (\ref{hpgde}) with
a dihedral monodromy group, say, the $\hpgo21$ solutions in (\ref{eq:dih12})--(\ref{eq:dih32}).
The functions $y_1^2$, $y_1y_2$, $y_2^2$ satisfy the same third order Fuchsian equation, 
called the {\em symmetric tensor square} of the second order equation under consideration.
The monodromy group of the symmetric square is the same dihedral group, 
but the corresponding monodromy representation is reducible \cite[Theorem 4.1]{singulm1}.
Hence an elementary solution of the symmetric square equation is expected.

Clausen's formula (\ref{eq:clausen}) identifies the tensor square equation of a 
rather general hypergeometric equation as a generalized hypergeometric equation
for $\hpgo32$ functions. It applies to the square of the dihedral $\hpgo21$ function
in (\ref{eq:diha}) with $k=0$, and to the square of 
\begin{equation} \label{eq:dih32a}
\hpg{2}{1}{\frac{a+1}{2}+k,\,\frac{1-a}{2}}{\frac{3}{2}+k}{\frac{z}{z-1}}.
\end{equation}
The latter function is related to the specialization $\ell=0$ of (\ref{eq:dih32}) by  
fractional-linear transformation (\ref{flinear3}). 
However, considering a similar transformation of (\ref{eq:dih12}) 
with $\ell=0$,  Clausen's identity with 
\begin{equation} \label{eq:dih32b}
\hpg21{\frac{a}2,-\frac{a}2-k}{\frac12-k}{\frac{z}{z-1}}^2
\end{equation}
leads to an apparently terminating $\hpgo32$ series with $A+B=-k$ in (\ref{eq:clausen}). 
Then Clausen's identity is false if the $\hpgo32$ series is interpreted as a terminating sum. 
Nevertheless, the terminating $\hpgo32$ sum still satisfies the same symmetric square equation;
it gives the expectable elementary solution. 

Here are correct identities that relate the terminating $\hpgo32$ series to squares of dihedral
$\hpgo21$ functions considered here, with $\ell=0$ or $k=0$. The $\hpgo32$ sum in (\ref{eq:symsq2})
is a rational function rather than a series around $z=1$.
\begin{lemma} \label{th:symsq}
For any $a\in\CC\setminus\ZZ$ and any integer $k\ge 0$ or $\ell\ge 0$ we have the following identities
in a neighbourhood of $z=0$: 
\begin{eqnarray} \label{eq:symsq}
&&\hspace{-30pt}(1-z)^{-a}\,\hpg32{-k,\,a,-a-2k}{-2k,\,\frac12-k}{\frac{z}{z-1}}=\nonumber\\
&&\hpg{2}{1}{\frac{a}{2},\frac{a+1}{2}}{\frac{1}{2}-k}{\,z}^2
-{ \frac{2^{4k}\,k!^4\,(a)_{2k+1}^{2}}{(2k)!^2(2k+1)!^2} }\, 
z^{2k+1}\,\hpg{2}{1}{\frac{a+1}{2}+k,\,\frac{a}{2}+k+1}{\frac{3}{2}+k}{\,z}^2,\\ \label{eq:symsq2}
&&\hspace{-30pt}\frac{(2\ell)!^2}{2^{4\ell}\;\ell!^2}\,(1-z)^{-a}
\hpg32{-\ell,\;a,-a-2\ell}{-2\ell,\,\frac12-\ell}{\frac1{1-z}}=\nonumber \\
&&{\textstyle \left(\frac{a+1}2\right)_\ell^2}\;\hpg21{\frac{a}2,\frac{a+1}2+\ell}{\frac12}{\,z}^2
- {\textstyle 4\left(\frac{a}2\right)_{\ell+1}^{2}}\,z\,\hpg21{\frac{a+1}2,\frac{a}2+\ell+1}{\frac32}{\,z}^2.
\end{eqnarray}
\end{lemma}
\begin{proof} We evaluate (\ref{eq:dih32b}) by taking a limit in Clausen's identity:
\begin{eqnarray*} \label{eq:exclausen}
\hpg21{A,\;-A-k}{\frac12-k}{\,x}^2 \hspace{-6pt} \equal
\lim_{B\to -A-k} \hpg21{A,\;B}{A+B+\frac12}{\,x}^2,\\
\equal \lim_{B\to -A-k} \hpg32{2A,\,2B,\,A+B}{2A+2B,A+B+\frac12}{\,x}.
\end{eqnarray*}
The latter limit is equal to
\begin{eqnarray*}
\hpg32{2A,\,-2A-2k,\,-k}{\frac12-k,\, -2k}{\,x}+ \hspace{240pt}\\
\frac{(2A)_{2k+1}(-2A\!-\!2k)_{2k+1}(-k)_k\,k!}{(\frac12-k)_{2k+1}\,(-2k)_{2k}\cdot 2}
\frac{x^{2k+1}}{(2k+1)!}\,\hpg32{2A+2k+1,1-2A,k+1}{\frac32+k,\,2k+2}{x}.
\end{eqnarray*}
Here the first $\hpgo32$ series is interpreted as terminating.
We apply Clausen's identity to the second $\hpgo32$ series,
and collect the $\hpgo21(x)^2$ terms on one side:
\begin{eqnarray} \label{eq:pre3f2}
\hpg32{2A,\,-2A-2k,\,-k}{\frac12-k,\, -2k}{\,x}=\hpg21{A,\;-A-k}{\frac12-k}{\,x}^2+
\hspace{100pt} \nonumber \\ 
\frac{\left(2A\right)_{2k+1}^2\,k!^2\,x^{2k+1}}
{2\left(\frac12\right)_k\!\left(\frac12\right)_{k+1}\!(2k)!\,(2k+1)!}\,
\hpg21{A+k+\frac12,\frac12-A}{k+\frac32}{x}^2.
\end{eqnarray}
Now we substitute $A=a/2$, $x=z/(z-1)$, use (\ref{eq:dfact}) to eliminate
a couple of Pochhammer symbols, 
and apply Pfaff's transformation (\ref{flinear3}) 
to the $\hpgo21$ series to get (\ref{eq:symsq}).

Formula (\ref{eq:symsq2}) can be proved by transforming the $\hpgo21$ functions
in (\ref{eq:pre3f2}) to 
\begin{equation*}
\hpg21{A,\;-A-k}{\frac12}{1-x}, \qquad \sqrt{1-x\,}\,\hpg21{A+\frac12,\frac12-A-k}{\frac32}{1-x}
\end{equation*}
using connection formula \cite[(2.3.13)]{specfaar} and Euler's transformation (\ref{flinear1}),
then applying the substitutions \mbox{$A=a/2$}, $k=\ell$, $x=1/(1-z)$ and Pfaff's transformation (\ref{flinear3}).
\end{proof}

The formulas of Lemma \ref{th:symsq} 
can be amalgamated to a general formula without the restriction $k\ell=0$.
From the results in \cite{gclausen} it follows that the elementary solution
(of the symmetric square equation) can be generally expressed using a terminating
$F^{2:1;1}_{1:1;1}$ sum defined in (\ref{eq:gap211}). 
\begin{theorem} \label{th:genelem}
For any $a\in\CC\setminus\ZZ$ and non-negative integers $k,\ell$, the following identity
holds in a neighbourhood of $z=0$: 
\begin{eqnarray} \label{eq:genelem}
 &&\hspace{-32pt}(1-z)^{-a}\,
F^{2:1;1}_{1:1;1}\!\left(\left. {a; -a-2k-2\ell; -k, -\ell\atop \frac12-k-\ell; \, -2k, -2\ell} \right| 
\frac{z}{z-1},\frac1{1-z}\right)=\nonumber \\
&&\frac{\left(\frac{a+1}2\right)_\ell\left(\frac{a+1}2+k\right)_\ell}
{\left(\frac12\right)_\ell\left(k+\frac12\right)_\ell}\;\hpg21{\frac{a}2,\frac{a+1}2+\ell}{\frac12-k}{\,z}^2-
\nonumber \\
&&\frac{\left(\frac{a+1}2\right)_k\left(\frac{a+1}2+\ell\right)_k\left(\frac{a}2\right)^2_{k+\ell+1}}
{\left(\frac12\right)_k\left(\frac12\right)_{k+1}^2\left(\frac12\right)_\ell\left(\frac12\right)_{k+\ell}}
\,z^{2k+1}\,\hpg21{\frac{a+1}2+k,\frac{a}2+k+\ell+1}{\frac32+k}{\,z}^2.
\end{eqnarray}
\end{theorem}
\begin{proof}
This is a reformulation of \cite[(11)]{gclausen}. The key result in \cite{gclausen} is that 
the univariate functions
\begin{eqnarray*}  \label{eq:spap0} 
\hpg21{a,\;b}{c}{\,x}^2, \qquad F^{2:1;1}_{1:1;1}\!\left(\left. {2a;\,2b;\;c-\frac12,\,a+b-c+\frac12 
\atop a+b+\frac12;\, 2c-1, 2a+2b-2c+1} \right|  x,1-x\right), \nonumber\\
\label{eq:spap1} x^{1-c}(1\!-\!x)^{c-a-b-\frac12}
F^{1:2;2}_{2:0;0}\!\left(\left. {\!\frac12; a-b+\frac12,a+b-c+\frac12; 
b-a+\frac12, c-a-b+\frac12 \atop c;\,2-c} 
\right|  x,\frac{x}{x\!-\!1}\right)
\end{eqnarray*}
satisfy the same Fuchsian equation of order 3.
From here the following direct generalization of Clausen's identity is derived \cite[(10)]{gclausen}:
\begin{eqnarray} \label{ge:clausen}
\hpg21{a,\;b}{a+b+\ell+\frac12}{\,x}^2= \hspace{215pt} \nonumber\\ 
\frac{(\frac12)_\ell\,(a+b+\frac12)_\ell}{(a+\frac12)_\ell\,(b+\frac12)_\ell}\,
F^{2:1;1}_{1:1;1}\!\left(\left. {2a;\,2b;\;a+b+\ell,\,-\ell 
\atop a+b+\frac12;\, 2a+2b+2\ell, -2\ell} \right|  x,1-x\right).
\end{eqnarray}
Identity (\ref{eq:genelem}) follows from (\ref{ge:clausen}) by taking the limit
$b\to-a-k-l$, similarly as in the proof of Lemma \ref{th:symsq}. 
\end{proof}

Besides, \cite{gclausen} proves
\begin{eqnarray} \label{eq:altaltclau}
&&\hspace{-32pt} \hpg{2}{1}{a,\,b\,}{c}{\,x\,}\,\hpg{2}{1}{1+a-c,\,1+b-c}{2-c}{\,x}=
(1-x)^{c-a-b-\frac12}\times\nonumber\\
&&F^{1:2;2}_{2:0;0}\!\left(\left. {\frac12; a-b+\frac12,a+b-c+\frac12; 
b-a+\frac12, c-a-b+\frac12 \atop c;\,2-c} 
\right|  x,\,\frac{x}{x-1}\right).
\end{eqnarray}
The case $c=a+b+\frac12$ of this formula is a well-known \cite[pg.~116]{specfaar} companion
to Clausen's identity,  due to Chaundy \cite{chaundy58}. 
Up to a power factor, the functions in (\ref{eq:altaltclau}) satisfy the same tensor square equation
as  $\displaystyle\hpg21{a,b}{c}{z}^2$.
Formula (\ref{eq:altaltclau}) can be applied with terminating $F^{1:2;2}_{2:0;0}$ series 
without any ambiguity. In particular, within our context we have 
\begin{eqnarray} \label{ge:altclaust}
\hpg{2}{1}{\frac{a}2,\,\frac{a+1}2+\ell}{a+k+\ell+1}{1-z}
\hpg{2}{1}{-\frac{a}2-k-\ell,\,\frac{1-a}2-k}{1-a-k-\ell}{1-z}=
\hspace{58pt}  \nonumber\\
z^{k}\,F^{1:2;2}_{2:0;0}\!\left(\left. {\frac12;\; k+1, \ell+1;\,-k,-\ell
\atop a+k+\ell+1;1-a-k-\ell} 
\right|  1-\frac1z,1-z\right).
\end{eqnarray}
This is a polynomial in $z$, the same as $f _{s,t}(1-z)$ in \cite[Section 3]{yoshidih} with 
\[
(s,t)=\left(\frac{k+\ell}2,\frac{k-\ell}2\right) \mbox{ or }
\left(\frac{k-\ell-1}2,\frac{k+\ell+1}2\right), \quad 
\mbox{and} \quad \alpha=a+k+\ell.
\]
Theorem 3.3 in \cite{yoshidih} states that the degree of the polynomial is $k+\ell$.

The $F^{1:2;2}_{2:0;0}$ sum in (\ref{ge:altclaust}) 
simplifies to a terminating $\hpgo32$ sum if $k=0$ or $\ell=0$. For example,
\begin{eqnarray} \label{eq:altaltclau2}
\hpg{2}{1}{\!\frac{a}{2},\frac{a+1}{2}+\ell}{a+\ell+1}{1\!-\!z}
\hpg{2}{1}{\!-\frac{a}{2}-\ell,\frac{1-a}{2}}{1-a-\ell}{1\!-\!z}=
\hpg32{\frac12,\,-\ell,\,\ell+1}{\!1-a-\ell,1+a+\ell}{1\!-\!z}. \quad
\end{eqnarray}
The terminating $F^{2:1;1}_{1:1;1}$ and  $F^{1:2;2}_{2:0;0}$ sums in 
(\ref{eq:genelem}) and (\ref{ge:altclaust}) are related by formula (\ref{eq:reverse32}).
Up to a constant multiple, the two sums are the same but the summation order is reversed
in both directions. This is consistent with one-dimensionality of the space of elementary solutions
for the symmetric tensor square equation. 

Technically, the elementary $F^{2:1;1}_{1:1;1}$ or  $F^{1:2;2}_{2:0;0}$ solutions are
values of degree 2 semi-invariants of the dihedral monodromy group 
(or the differential Galois group) of starting hypergeometric equation \cite{singulm1}. 
Formulas (\ref{eq:genelem}) and (\ref{ge:altclaust}) identify the semi-invariants
as quadratic expressions in $\hpgo21$ solutions. 
There are higher degree semi-invariants if the monodromy group is a finite dihedral group;
see Remark \ref{rm:minvars} below. 


\section{Transformations of dihedral functions} 
\label{sec:trihedr}

In the following two sections we characterize algebraic transformations between 
dihedral hypergeometric functions. In the classification \cite{algtgauss} of all 
transformations between Gauss hypergeometric functions, the following pull-back transformation
of general dihedral hypergeometric equations comes up:
\begin{equation} \label{eq:ledpb} \textstyle
E\left(\frac12,\frac12,a\right)\stackrel{n}{\longleftarrow} E\left(\frac12,\frac12,na\right).
\end{equation}
These transformations have any degree $n$, and have a free parameter $a$. 
There are more transformations if  $a\in\QQ\setminus\ZZ$; we consider them in Section \ref{adihedral}.

In the notation of \cite{algtgauss},
the branching pattern for a pull-back covering $\varphi:\PP^1_x\to\PP^1_z$ 
of transformation (\ref{eq:ledpb}) must be
\begin{eqnarray*}
1+2+2+\ldots+2=n=1+2+2+\ldots+2, && \mbox{if $n$ is odd},\\
1+1+2+2+\ldots+2=n=2+2+\ldots+2,\hspace{17pt} && \mbox{if $n$ is even}.
\end{eqnarray*}
This means that a fiber (above a singular point with the local exponent difference $a$) 
has a single point with the branching order $n$,
and that there are $n-1$ simple branching points in the other two fibers above singular points.
The covering is \mbox{$z=\varphi(x)$}, where $\varphi(x)$ is a rational function of degree $n$.
We assume that the local exponent  differences $1/2$ are assigned to the points $z=0$, $z=\infty$ 
and $x=0$, $x=\infty$. We assume that the point $x=0$ lies above the point $z=0$, 
so that $x$ divides the numerator of $\varphi(x)$. 
An expected hypergeometric identity is
\begin{equation}
\hpg{2}{1}{\frac{an}2,\,\frac{an+1}2}{1/2}{\,x}=\theta(x)\,
\hpg{2}{1}{\frac{a}2,\,\frac{a+1}2}{1/2}{\varphi(x)}.
\end{equation}
Here we present the explicit covering and corresponding transformations.
The two branching patterns are asserted by polynomial identity (\ref{eq:dhtrb}).
\begin{theorem} \label{th:transf}
For a positive integer $n$, let us define the polynomials
\begin{equation} \label{thetasd}
\theta_1(x)=\sum_{k=0}^{\lfloor n/2 \rfloor} {n\choose 2k}x^k,
\qquad \theta_2(x)=\sum_{k=0}^{\lfloor (n-1)/2 \rfloor}
 {n\choose 2k+1} x^k.
\end{equation}
Then the rational function $\varphi(x)=x\,\theta_2(x)^2\big/\theta_1(x)^2$ 
realizes pull-back transformation $(\ref{eq:ledpb})$,
and we have the following identities
\begin{eqnarray} 
\label{eq:dhtra} \left(1-x\right)^n \equal \theta_1\!\left(x^2\right)-x\,\theta_2\!\left(x^2\right),\\
\label{eq:dhtrb} \left(1-x\right)^n \equal \theta_1(x)^2-x\,\theta_2(x)^2,\\
\label{eq:dhtr1} \hpg{2}{1}{\frac{na}{2},\,\frac{na+1}{2}}{\frac{1}{2}}{\,x}
\equal\theta_1(x)^{-a}
\,\hpg{2}{1}{\frac{a}{2},\,\frac{a+1}{2}}{\frac{1}{2}}
{\frac{x\,\theta_2(x)^2}{\theta_1(x)^2}},\\
\label{eq:dhtr2}  \hpg{2}{1}{\frac{na+1}{2},\,\frac{na}{2}+1}{\frac32}{\,x}\equal
\theta_1(x)^{-a-1}\,\frac{\theta_2(x)}{n}\,\hpg{2}{1}{\frac{a+1}{2},\,
\frac{a}{2}+1}{\frac32}{\frac{x\,\theta_2(x)^2}{\theta_1(x)^2}},\\
\label{eq:dhtr3}  \hpg{2}{1}{\frac{na}{2},\,-\frac{na}{2}}{\frac12}{\frac{x}{x-1}}\equal
\hpg{2}{1}{\frac{a}{2},\,-\frac{a}{2}}{\frac12}{-\frac{x\,\theta_2(x)^2}{(1-x)^n}},\\
\label{eq:dhtr4} \hpg{2}{1}{\frac{na}{2},\frac{na+1}{2}}{na+1}
{1-x}\equal\left(\frac{\theta_1(x)}{2^{n-1}}\right)^{-a}
\hpg{2}{1}{\frac{a}{2},\,\frac{a+1}{2}}{a+1}
{\frac{(1-x)^n}{\theta_1(x)^2}}.
\end{eqnarray} 
\end{theorem}
\begin{proof} The polynomials $\theta_1(x)$, $\theta_2(x)$ satisfy
\begin{equation} \label{eq:sqrtheta}
\left(1-\sqrt{x}\right)^n=\theta_1(x)-\sqrt{x}\,\theta_2(x).
\end{equation}
After the substitution $x\mapsto x^2$ we get (\ref{eq:dhtra}).
Multiplication of (\ref{eq:sqrtheta}) with
its own conjugate $\sqrt{x}\mapsto-\sqrt{x}$ version 
gives (\ref{eq:dhtrb}).

Formula (\ref{eq:dhtr1}) is easily obtainable using (\ref{dihedr2}):
\begin{eqnarray*}
\hpg{2}{1}{\frac{na}{2},\,\frac{na+1}{2}}{\frac{1}{2}}{\,x} \equal
\frac{(1-\sqrt{x})^{-na}+(1+\sqrt{x})^{-na}}{2} \nonumber\\
\equal \frac{\big(\theta_1(x)-\sqrt{x}\,\theta_2(x)\big)^{-a}
+\big(\theta_1(x)+\sqrt{x}\,\theta_2(z)\big)^{-a}}{2}\\
\equal \theta_1(x)^{-a}
\,\hpg{2}{1}{\frac{a}{2},\,\frac{a+1}{2}}{\frac{1}{2}}{\frac{x\,\theta_2(x)^2}{\theta_1(x)^2}}.
\end{eqnarray*}
By part 2 of \cite[Lemma 2.1]{algtgauss}, the covering 
$z=\varphi(x)=x\,\theta_2(x)^2\big/\theta_1(x)^2$
gives a pull-back (\ref{eq:ledpb}) between the corresponding differential equations.
Formula (\ref{eq:dhtr2}) follows from  \cite[Lemma 2.3]{algtgauss} applied to 
(\ref{eq:dhtr1}). The last two formulas follow from a standard identification of
local hypergeometric solutions (at $x=0$ and $x=1$, respectively) related by a pull-back
transformation.
\end{proof}

The polynomials $\theta_1(x)$, $\theta_2(x)$ in  (\ref{thetasd}) can be written as 
$\hpgo21$ hypergeometric sums:
\begin{equation}
\theta_1(x)=\hpg21{-\frac{n}2,-\frac{n-1}2}{1/2}{x},\qquad
\theta_2(x)=n\,\hpg21{-\frac{n-1}2,-\frac{n-2}2}{3/2}{x}.
\end{equation}
The peculiarly similar identities (\ref{eq:dhtra}) and (\ref{eq:dhtrb}) can be written as follows:
\begin{eqnarray}
(1-x)^n\equal \hpg21{-\frac{n}2,-\frac{n-1}2}{1/2}{x^2}
-n\,x\,\hpg21{-\frac{n-1}2,-\frac{n-2}2}{3/2}{x^2},\\ 
(1-x)^n  \equal \hpg21{-\frac{n}2,-\frac{n-1}2}{1/2}{x}^2
-n^2x\,\hpg21{-\frac{n-1}2,-\frac{n-2}2}{3/2}{x}^2.
\end{eqnarray}
The latter identity is the special case $a=-n$, $k=0$ of Lemma \ref{th:symsq}.
The polynomials $\theta_1(x)$, $\theta_2(x)$ are related to the Tchebyshev polynomials 
(of the first and the second kind):
\begin{equation} \label{eq:tchebysh}
T_n(x)=x^n\,\theta_1\!\left(\frac{x^2-1}{x^2}\right),\qquad
U_{n-1}(x)=x^{n-1}\,\theta_2\!\left(\frac{x^2-1}{x^2}\right).
\end{equation}

Formula (\ref{eq:dhtr3}) can be written, for odd $n$, as
\begin{eqnarray}
\hpg{2}{1}{\frac{na}{2},-\frac{na}{2}}{\frac12}{\,x}\equal
\hpg{2}{1}{\frac{a}{2},-\frac{a}{2}}{\frac12}{n^2x\,\hpg21{\frac{1-n}2,\frac{1+n}2}{3/2}{x}^{\!2}},
\end{eqnarray}
and for even $n$ as 
\begin{eqnarray}
\hpg{2}{1}{\frac{na}{2},-\frac{na}{2}}{\frac12}{\,x}\equal
\hpg{2}{1}{\frac{a}{2},-\frac{a}{2}}{\frac12}{n^2x(1-x)\hpg21{1-\frac{n}2,1+\frac{n}2}{3/2}{x}^{\!2}}.
\end{eqnarray}

The transformation of Theorem \ref{th:transf} is unique with the prescribed ramification pattern
(with $x=0$ lies above $z=0$ where the local exponent difference is $1/2$), 
because such an transformation with $\varphi$ normalized as described  must identify
the explicit solutions \mbox{$(1-\sqrt{x})^{-na}$} and \mbox{$(1-\sqrt{z})^{-a}$}, 
and there is only one way to identify them.

Transformation (\ref{eq:ledpb}) can be composed with
quadratic transformation (\ref{eq:quatriv}) with $k=0$, to a transformation
\begin{equation} \label{eq:ledpc}
\textstyle
E\left(\frac12,\frac12,a\right)\stackrel{2n}{\longleftarrow}E(1,na,na).
\end{equation} 
The same transformation can be obtained by composing the same quadratic transformation and 
$E(1,a,a)\stackrel{n}{\longleftarrow}E(1,na,na)$ of \cite[Section 5]{algtgauss}. 
By the classification in \cite[Section 5]{algtgauss}, transformations (\ref{eq:ledpb}), (\ref{eq:ledpc})
and the quadratic transformations exhaust  all pull-back transformations of hypergeometric equations 
with an infinite dihedral monodromy group. 
There are higher degree pull-back transformations of dihedral hypergeometric functions only if
$a$ is a rational number; the monodromy group is then a finite dihedral group.
Klein's theorem \cite{klein77} states that any second order Fuchsian equation
with a finite monodromy group is a pull-back of one of a few {\em standard hypergeometric equations}
with the same monodromy group. This implies that any algebraic dihedral function
is a pull-back transformation of a solution of $E(1/2,1/2,1/m)$,
where $m$ is a positive integer. The relevant Klein pull-back
transformations are
\begin{equation} \label{eq:kleinled} \textstyle
E\left(\frac12,\frac12,\frac1m\right)\stackrel{d}{\longleftarrow}E\left(k+\frac12,\ell+\frac12,\cc\right),
\end{equation}
where $\cc$ is a rational number whose denominator is $m$. 
We consider these transformations in the following section.
The covering degree turns out to be $d=(k+\ell)m+\lambda m$.

\section{Transformations of algebraic dihedral functions}
\label{adihedral}

If $a$ is a rational number but not an integer, the dihedral Gauss hypergeometric functions
are algebraic functions. There are more pull-back transformations (\ref{hpgtransf}) 
between these functions. In particular, there are transformations (\ref{eq:kleinled})
implied by celebrated Klein's theorem. Here we consider Klein's pull-back transformations
from $E(1/2,1/2,1/m)$, where $m$ is a positive integer greater than $1$. 
The case $m=1$ gives degenerate or logarithmic dihedral hypergeometric functions;
their transformations are commented in Remark \ref{rm:degtrans} below.

Klein's pull-back transformations for second order linear equations with a finite monodromy group
can be found by computing values of monodromy semi-invariants as solutions of appropriate symmetric tensor powers of the given equation \cite{kleinvhw}. 
For Gauss hypergeometric functions with a finite monodromy group, Klein's pull-back transformations
can be found by using contiguous relations and a data base of simplest explicit expressions of Gauss
hypergeometric functions of each Schwarz type \cite{talggaus}. The approach in \cite{talggaus} 
simplifies greatly for algebraic dihedral hypergeometric functions, as we demonstrate here. In particular, we do not have to use contiguous relations, since we already have general explicit expressions for any dihedral hypergeometric functions.

\begin{theorem} \label{th:genftrans}
Let $H_1$ denote hypergeometric equation $E(k+1/2,\ell+1/2,n/m)$, 
where $k,\ell,m,n$ are positive integers, $m>1$ and $\gcd(n,m)=1$.
The monodromy group of $H_1$ is the finite dihedral group with $2m$ elements.
Let us denote
\begin{equation} \label{eq:gdef1}
G(\sqrt{x})= 
x^{k/2}\,\app3{k+1,\ell+1; -k,-\ell}{1+\frac{n}m}{\frac{\sqrt{x}+1}{2\sqrt{x}},\frac{1+\sqrt{x}}2}.
\end{equation}
This is a polynomial in $\sqrt{x}$. We can write 
\begin{equation} \label{eq:gmpower}
\left(1+\sqrt{x}\right)^{n}G(\sqrt{x})^m=\Theta_1(x)+x^{k+\frac12}\,\Theta_2(x),
\end{equation}
so that $\Theta_1(x)$ and $\Theta_2(x)$ are polynomials in $x$. 
Klein's pull-back covering for $H_1$ is given by the rational function 
$\Phi(x)=x^{2k+1}\,\Theta_2(x)^2/\Theta_1(x)^2$.
The covering is unique up to M\"obius transformations,
and its degree is equal to $(k+\ell)m+n$.
\end{theorem}
\begin{proof}  
The monodromy group is clarified by quadratic transformation (\ref{eq:qled}).
The corresponding standard hypergeometric equation $H_0$ is $E(1/2,1/2,1/m)$. 
The branching pattern above the point with the local exponent difference
$1/m$ must be $n+m+m+\ldots+m$. Above the two points with the local exponent 
difference $1/2$ there must
be single points of ramification order $2k+1$ and $2\ell+1$ (which may lie in the same fiber or not), 
and the remaining points are simple ramification points. The degree of Klein's pull-back covering
for $H_1$ follows from Hurwitz formula, or immediately from part 2 of \cite[Lemma 2.5]{algtgauss}:
\begin{equation}
d=\frac{(k+\frac12)+(\ell+\frac12)+\frac{n}m-1}{\frac12+\frac12+\frac1m-1}=(k+\ell)\,m+n.
\end{equation} 

As explained in \cite[Section 2]{talggaus}, Klein's pull-back covering can be seen as the composition
$s_0^{-1}\circ s$, where $s_0$ is a Schwarz map for the standard equation $H_0$, and $s$ is the corresponding Schwarz map for the given equation $H_1$. The inverse Schwarz map $s_0^{-1}$
for the standard equation is a rational function (of degree $2m$), and it cancels the monodromy of $s$. 
Recall that a {\em Schwarz map} for a second order linear differential equation is the ratio of any
two distinct solutions of it. The image of the complex upper half-plane under $s_0$ is a spherical triangle 
with the angles $(\pi/2,\pi/2,\pi/m)$.

Slightly differently from \cite[Algorithm 3.1]{talggaus}, let $z$, $x$ denote the variables for $H_0$
and $H_1$, respectively. 
By M\"obius transformations we can choose three points on any projective line freely.
We assign the local exponent difference $1/2$ to the points $z=0$ and $z=\infty$,
and the local exponent difference $1/m$ to $z=1$. On the $x$-projective line, we assign
$k+1/2$, $\ell+1/2$, $n/m$ to the points $x=0$, $x=\infty$ and $x=1$, respectively. 
We choose the point $x=0$ to lie above $z=0$; the point $x=1$ must lie above $z=1$; 
the point $x=\infty$ lies above $z=0$ or $z=\infty$.
The Darboux pull-back coverings (as defined in \cite{dalggaus} and \cite{talggaus}) for both equations are simply the quadratic transformations that reduce the dihedral monodromy group to $\ZZ/m\ZZ$. 
The Darboux coverings are defined by the functions $\sqrt{z}$ and $\sqrt{x}$.

We have the following two solutions of the standard equation $H_0$:
\begin{eqnarray} \label{sdihedr1}
\hpg{2}{1}{-\frac{1}{2m},\,\frac{1}{2}-\frac{1}{2m}}{1-\frac{1}{m}}
{1-z} & = & \left(\frac{1+\sqrt{z}}{2}\right)^{1/m}.\\
\label{sdihedr2}
(1-z)^{1/m}\;\hpg{2}{1}{\frac{1}{2m},\,\frac{1}{2}+\frac{1}{2m}}{1+\frac{1}{m}}
{1-z} & = & \big(2-2\sqrt{z}\big)^{1/m}.
\end{eqnarray}
We choose the Schwartz map to be 
\begin{equation}
s_0=\left(\frac{1+\sqrt{z}}{1-\sqrt{z}}\right)^{1/m},
\end{equation}
which is a quotient of the two solutions up to the constant multiple $4^{1/m}$. 
The inverse Schwartz map is 
\begin{equation} \label{eq:invsch}
z=\frac{(s_0^m-1)^2}{(s_0^m+1)^2}. 
\end{equation}

The hypergeometric solutions of $H_1$ corresponding to (\ref{sdihedr1})--(\ref{sdihedr2}) are
(\ref{eq:diha}) and 
\begin{equation}
(1-z)^{-a-k-\ell}\,\hpg{2}{1}{-\frac{a}{2}-k-\ell,\,\frac{1-a}{2}-k}{1-a-k-\ell}{1-z},
\end{equation}
with $a=-n/m-k-\ell$ and $z=x$. The correspondence is up to the same
fractional power factor and possibly different scalar factor. 
The Schwartz map $s$ is equal, up to a constant factor, to
\begin{eqnarray*}
\quad\hpg{2}{1}{\!-\frac{n}{2m}-\frac{k+\ell}2,\frac{1-k+\ell}{2}-\frac{n}{2m}}{1-\frac{n}{m}}{1-x}
\left/ (1-x)^{n/m}\,
\hpg{2}{1}{\!\frac{n}{2m}-\frac{k+\ell}2,\frac{1-k+\ell}{2}+\frac{n}{2m}}{1+\frac{n}{m}}{1-x}
\right..
\end{eqnarray*}
After applying (\ref{eq:diha}) and (\ref{eq:f3f3rel})  
to the numerator $\hpgo21$ function
and just (\ref{eq:diha}) to the denominator $\hpgo21$ function we conclude that 
the Schwartz map $s$ is equal,  up to a constant factor, to
$\widetilde{G}\left(\sqrt{x}\right)/\widetilde{G}\left(-\sqrt{x}\right)$, where
\begin{equation} \label{kleing0}
\widetilde{G}(t)=(1+t)^{n/m}\,\app3{k+1,\,\ell+1;\, -k,\,-\ell}{1+n/m}{\frac{t+1}{2\,t},\frac{1+t}2}.
\end{equation}
The constant multiple is restricted by compatibility of the conjugating automorphisms 
\mbox{$\sqrt{z}\mapsto-\sqrt{z}$} and $\sqrt{x}\mapsto-\sqrt{x}$ of the two Darboux coverings. 
This automorphism acts on $s_0$ as $s_0\mapsto1/s_0$; the action on $s$ must be the same.
It follows that $s$ is equal to $\widetilde{G}\left(\sqrt{x}\right)/\widetilde{G}\left(-\sqrt{x}\right)$,
possibly up to multiplication by $-1$.
From explicit form (\ref{eq:invsch}) of $s_0^{-1}$ we conclude that
\begin{equation} \label{eq:kleing2}
z=\left( \frac{\widetilde{G}\left(\sqrt{x}\right)^m-\widetilde{G}\left(-\sqrt{x}\right)^m}
{\widetilde{G}\left(\sqrt{x}\right)^m+\widetilde{G}\left(-\sqrt{x}\right)^m} \right)^{\pm2}.
\end{equation}
Since we chose $x=0$ to lie above $z=0$, the correct exponent is $2$ rather than $-2$.
Uniqueness already follows.
We indeed get a rational function in $x$ on the right-hand side of (\ref{eq:kleing2}); 
the function is odd with respect to $x\mapsto-x$. The left-hand side of (\ref{eq:gmpower}) 
is equal to $\widetilde{G}\left(\sqrt{x}\right)^m$.
The difference $\widetilde{G}\left(\sqrt{x}\right)^m-\widetilde{G}\left(-\sqrt{x}\right)^m$ 
must be of order $k+1/2$, since this is the local exponent difference at $x=0$. 
\end{proof}

\begin{corollary} \label{th:genftrans2}
We can obtain the same rational function $\Phi(x)$ of the same form 
if we replace the left-hand side of $(\ref{eq:gmpower})$ by
\begin{equation}  \label{eq:gdef2}
\left(1+\sqrt{x}\right)^{n+(k+\ell)m}\,\app2{-\frac{n}m-k-\ell;-k,-\ell}{-2k,-2\ell}{\frac{2\sqrt{x}}{1+\sqrt{x}},\frac2{1+\sqrt{x}}}^m
\end{equation}
and define  $\Theta_1(x)$, $\Theta_2(x)$ in the same way.
\end{corollary}
\begin{proof} Appell's $F_2$ and $F_3$ functions in (\ref{eq:gdef1}) and (\ref{eq:gdef2})
are related by transformation (\ref{eq:r2f3rel}), hence the two versions of 
$\Theta_1(x)$ and $\Theta_2(x)$ differ by the same constant. It is instructive to observe 
that the same Schwarz map $s_0$ can be defined in terms of solutions (\ref{dihedr2}) and
the similar solution \cite[(1.3)]{tdihedral} 
of $H_0$, and then consider possible correspondence to $H_1$-solutions (\ref{eq:dih12}) and (\ref{eq:dih32}).\end{proof}

The following identity 
is an analogue of (\ref{eq:dhtrb}).
\begin{corollary} \label{th:invsrel}
We have the identity
\begin{equation}  \label{eq:insrel}
\Theta_1(x)^2-x^{2k+1}\,\Theta_2(x)^2=C\,(1-x)^n\,\Psi(x)^m,
\end{equation}
where
\begin{eqnarray} \label{eq:psii1}
\!\! \Psi(x)\equal  \frac{(-1)^\ell\left(1+\frac{n}m\right)_{k+\ell}\left(1-\frac{n}m\right)_{k+\ell}
\left(\frac{1-k-\ell}2-\frac{n}{2m}\right)_\ell}{4^{k+\ell}\left(\frac{1+k-\ell}2-\frac{n}{2m}\right)_\ell}
\times\nonumber\\ &&x^k\, F^{1:2;2}_{2:0;0}\!\left(\left. {\frac12;\; k+1, \ell+1;\,-k,-\ell}
\atop {1+\frac{n}m;1-\frac{n}m} \right| 1-\frac1x,1-x \right) \\
\equal \frac{\left(\frac12\right)_{k+\ell}\left(\frac12\right)_k\left(\frac12\right)_\ell
\left(\frac{1-k-\ell}2-\frac{n}{2m}\right)_\ell}{\left(\frac{1+k-\ell}2-\frac{n}{2m}\right)_\ell} 
\times\nonumber\\ && (1-x)^{k+\ell}\,
F^{2:1;1}_{1:1;1}\!\left(\left. {\!\frac{n}m-k-\ell; -\frac{n}m-k-\ell; -k, -\ell\atop \frac12-k-\ell;\,-2k, -2\ell}
\right| \frac{x}{x-1},\frac1{1-x}\right), \quad
\end{eqnarray}
and
\begin{eqnarray*}
C=4^{n+m(k+\ell)}\left/{\textstyle \left(1+\frac{n}m\right)_{k+\ell}^{2m} } \right.
&& \mbox{in the setting of Theorem } \ref{th:genftrans},\\
C=4^n\left/{\textstyle \left(\frac12\right)_k^{2m}\left(\frac12\right)_\ell^{2m}} \right. \hspace{28pt}
&& \mbox{in the setting of Corollary } \ref{th:genftrans2}.
\end{eqnarray*}
\end{corollary}
\begin{proof} Multiplication of both sides of $(\ref{eq:gmpower})$ respectively with the conjugate 
$\sqrt{x}\mapsto -\sqrt{x}$ of itself gives
\begin{eqnarray*}
(1-x)^n\,G(\sqrt{x})\,G(-\sqrt{x})=\Theta_1(x)^2-x^{2k+1}\,\Theta_2(x)^2.
\end{eqnarray*}
To evaluate $G(\sqrt{x})\,G(-\sqrt{x})$, we use formulas (\ref{eq:gdef1}), (\ref{eq:diha}), (\ref{ge:altclaust}),
and optionally (\ref{eq:reverse32}), (\ref{eq:gdef2}), (\ref{eq:r2f3rel}). 
To check the constant $C$, one may use the explicit expressions in Remark \ref{rm:minvars} below.
\end{proof}

Here are some algebraic transformations of dihedral hypergeometric functions
implied by Theorem \ref{th:genftrans}, parallel to (\ref{eq:dhtr1})--(\ref{eq:dhtr4}):
\begin{eqnarray}
\hpg{2}{1}{\!-\frac{n}{2m}\!-\!\frac{k+\ell}2,\frac{1-k+\ell}2\!-\!\frac{n}{2m}}{\frac{1}{2}-k}{\,x} \!\equal\!\!
\left(\frac{\Theta_1(x)}{\Theta_1(0)}\right)^{1/m} \! 
\hpg{2}{1}{\!-\frac{1}{2m},\frac12\!-\!\frac{1}{2m}}{\frac{1}{2}}{\frac{x^{2k+1}\Theta_2(x)^2}{\Theta_1(x)^2}},
\qquad\\ 
\hpg{2}{1}{\!\frac{1+k-\ell}2\!-\!\frac{n}{2m},\frac{k+\ell}2\!+\!1\!-\!\frac{n}{2m}}{\frac{3}{2}+k}{\,x} \!\equal\!
\left(\frac{\Theta_1(x)}{\Theta_1(0)}\right)^{1/m-1}\,\frac{\Theta_2(x)}{\Theta_2(0)}\times \nonumber\\
&& \hpg{2}{1}{\frac12-\frac{1}{2m},\,1-\frac{1}{2m}}
{\frac{3}{2}}{\frac{x^{2k+1}\Theta_2(x)^2}{\Theta_1(x)^2}},\qquad\\
\hpg{2}{1}{\!-\frac{n}{2m}\!-\!\frac{k+\ell}2,\frac{n}{2m}\!-\!\frac{k+\ell}2}{\frac{1}{2}-k}{\frac{x}{x-1}} \!\equal\!
\sqrt{\frac{\Psi(x)}{\Psi(0)}}\,(1-x)^{-\frac12(k+\ell)}\times \nonumber\\ && 
\hpg{2}{1}{-\frac{1}{2m},\;\frac{1}{2m}}{\frac{1}{2}}{-\frac{x^{2k+1}\Theta_2(x)^2}{C(1\!-\!x)^{n}\Psi(x)^m}},\\
\hspace{-6pt} 
\hpg{2}{1}{\!-\frac{n}{2m}\!-\!\frac{k+\ell}2,\frac{1-k+\ell}2\!-\!\frac{n}{2m}}{1-\frac{n}{m}}{1\!-\!x} \!\equal\!
\left(\frac{\Theta_1(x)}{\Theta_1(1)}\right)^{1/m} \times \nonumber\\ && 
\hpg{2}{1}{\!-\frac{1}{2m},\frac12-\frac{1}{2m}}{1-\frac{1}{m}}{\frac{C(1\!-\!x)^{n}\Psi(x)^m}{\Theta_1(x)^2}}.
\end{eqnarray}

\begin{remark} \label{rm:schwarz} \rm
As in Theorem \ref{th:genftrans}, let $H_1$ denote a hypergeometric equation
with a finite dihedral monodromy group, and let $H_0$ denote a standard hypergeometric equation
with the same monodromy group. Klein's pull-back covering is a composition 
$s_0^{-1}\circ s$ of the Schwarz map $s:\PP^1_x\to\PP^1_s$ for $H_1$
with the corresponding inverse Schwarz map $s_0^{-1}$ for $H_0$.
In \cite[Section 4]{ochiay}, Ochiai and Yoshida consider the composition $s\circ s_0^{-1}$, 
which is a rational function (of $s_0$) as well. 
Following the construction of (\ref{eq:kleing2}), we can write this rational function
(denoted by $w$ in \cite{ochiay}) as
\begin{eqnarray}
\frac{\displaystyle
(-1)^ks_0^n\,\app3{k+1,\,\ell+1;\, -k,\,-\ell}{1+\frac{n}m}{\frac{s_0^m}{s_0^m-1},\frac{s_0^m}{s_0^m+1}}}
{\displaystyle\app3{k+1,\,\ell+1;\, -k,\,-\ell}{1+\frac{n}m}{\frac{1}{1-s_0^m},\frac{1}{1+s_0^m}}},
\end{eqnarray}
or equivalently,
\begin{eqnarray}
\frac{\displaystyle
s_0^{n+m(k+\ell)}\,\app2{-\frac{n}m-k-\ell;\, -k,\,-\ell}{-2k,-2\ell}{1-s_0^{-m},1+s_0^{-m}}}
{\displaystyle\app2{-\frac{n}m-k-\ell;\, -k,\,-\ell}{-2k,-2\ell}{1-s_0^m,1+s_0^m}}.
\end{eqnarray}
The denominator of this rational function is a constant multiple of the polynomial $Dih(s_0^m)$
of \cite{ochiay}.  Here is the identification of the notation of \cite{ochiay} and this paper
in the dihedral context:
\begin{eqnarray}
\left(\mu_0=\mu_1,\,\mu_\infty,\,k_0=k_1,\,n=k_\infty, p,\,\big[p_0,\,p_1,\,p_\infty\big],\,
\,q,\,r,\,z,\,\frac{f_1}{f_0} \,\right) \hspace{32pt} \nonumber \\
\mapsto \left(\, \frac12,\,\frac1m,\,2,\,m,\,k+\ell,\,\left[\ell,\,k,\,\frac{n-1}m\right]\!,\,
\frac{n}m,\,k-\ell,\,s_0,\,\sqrt{z} \, \right).
\end{eqnarray}
Following \cite[(3.6)]{ochiay} and formulas (\ref{eq:diha}), (\ref{eq:r2f3rel}), (\ref{eq:f2f2rel})
here, we indeed get
\begin{eqnarray}
Dih(t) \equal (1+t)^{\frac{n}m+k+\ell}\,
\hpg21{-\frac{n}{2m}-\frac{k+\ell}2,-\frac{n}{2m}+\frac{1-k+\ell}2}{1-\frac{n}m}{\frac{4t}{(1+t)^2}}\\
\equal \frac{(\frac{n}{2m}+\frac{1-k-\ell}2)_k\,(2k)!\,(2\ell)!}
{(\frac{n}{2m}+\frac{1-k+\ell}2)_k\,(1-\frac{n}{m})_{k+\ell}\,k!\,\ell!}\,
\app2{-\frac{n}{m}-k-\ell;\, -k,\,-\ell}{-2k,-2\ell}{1-t,1+t}.\qquad
\end{eqnarray}
Ochiai and Yoshida investigate the coefficients of $Dih(t)$.
In differential equation \cite[(4.6)]{ochiay} we recognize Heun equation (\ref{eq:heun}) 
with (\ref{heunhpg2b}).
Theorems 4.13 and 4.15 in \cite{ochiay} state the following expressions, with $\cc=-n/m$:
\begin{eqnarray*}
Dih(t) \!\equal\! \sum_{p=0}^{\lfloor\frac{k+\ell}2\rfloor} 
\frac{\big(\frac{\cc-k-\ell}2\big)_p\left(\frac{1-k-\ell}2\right)_p}
{\big(\frac{1+\cc}2\big)_{p}\;p!}\,
\hpgg43{-p,\,p+\frac{\cc-k-\ell}2,\frac{1+k-\ell}2,\frac{1-k+\ell}2}
{\frac12,\,1+\frac{\cc}{2},\,\frac{1-k-\ell}2} t^{2p}\\
&&\!+(\ell-k)\sum_{p=0}^{\lfloor\frac{k+\ell}2\rfloor} 
\frac{\big(\frac{\cc-k-\ell}2\big)_{p+1}\left(\frac{1-k-\ell}2\right)_p}{\big(\frac{1+\cc}{2}\big)_{p+1}\;p!}
\hpgg43{\!-p,p+1+\frac{\cc-k-\ell}2,\frac{1+k-\ell}2,\frac{1-k+\ell}2}
{\frac32,\,1+\frac{\cc}{2},\,\frac{1-k-\ell}2} t^{2p+1}\quad\\
\equal\! \sum_{p=0}^{\lfloor\frac{k+\ell}2\rfloor} 
\frac{\big(\frac{\cc-k-\ell}2\big)_p\,\left(-\frac{k+\ell}2\right)_p}{(1+\frac{\cc}{2})_{p}\;p!}\,
\hpgg43{-p,\,p+\frac{\cc-k-\ell}2,\,\frac{k-\ell}2,\,\frac{\ell-k}2}
{\frac12,\,\frac{1+\cc}{2},\,-\frac{k+\ell}2} t^{2p}+\\
&&\hspace{-4pt}(\ell-k)\sum_{p=0}^{\lfloor\frac{k+\ell}2\rfloor} 
\frac{\!\big(\frac{\cc-k-\ell}2\big)_{p+1}\left(1-\frac{k+\ell}2\right)_p}
{\frac{1+\cc}{2}\;\big(1+\frac{\cc}{2}\big)_{p}\;p!}
\hpgg43{\!-p,p+1+\frac{\cc-k-\ell}2,1+\frac{k-\ell}2,1+\frac{\ell-k}2}
{\frac32,\,\frac{3+\cc}{2},\,1-\frac{k+\ell}2} t^{2p+1}.
\end{eqnarray*}
The argument of the $\hpgo43(x)$ functions is $x=1$.
\end{remark}

\begin{remark} \rm \label{rm:minvars}
The polynomials $\Theta_1(x)$ and $\Theta_2(x)$ of Theorem \ref{th:genftrans} 
or Corollary \ref{th:genftrans2} appear as semi-invariants (of degree $m$) 
of the finite dihedral monodromy group. Let us consider the following solutions of
a dihedral hypergeometric equation, 
with $\cc=-n/m$:
\begin{eqnarray*}
W_1 \equal
\frac{2^{-\cc-k-\ell}\left(1+\cc\right)_{k+\ell}}{\left(\frac{1+\cc+k-\ell}2\right)_\ell}
\;\hpg21{\!\frac{\cc-k-\ell}2,\frac{1+\cc-k+\ell}{2}}{1+\cc}{1-x}\\
\equal {\textstyle \left(\frac12\right)_k}\;
\hpg21{\!\frac{\cc-k-\ell}2,\frac{1+\cc-k+\ell}{2}}{\frac12-k}{\,x}+ \\ & & 
\frac{(-1)^{k+1}\left(\frac{1+\cc-k-\ell}2\right)_k
\left(\frac{\cc-k-\ell}2\right)_{k+\ell+1}}
{\left(\frac{1+\cc-k-\ell}2\right)_\ell\;\left(\frac12\right)_{k+1}}\;
x^{k+\frac12}\hpg21{\!\frac{\cc-k-\ell}2,\frac{1+\cc-k+\ell}{2}}{\frac32+k}{\,x},\\
W_2 \equal \frac{(-1)^{\ell}\,2^{-\cc-k-\ell}\left(1-\cc\right)_{k+\ell}}
{\left(\frac{1+\cc-k-\ell}2\right)_\ell}\,(1-x)^{-\cc}
\;\hpg21{\frac{-\cc-k-\ell}2,\frac{1-\cc-k+\ell}{2}}{1-\cc}{1-x}\\
\equal {\textstyle \left(\frac12\right)_k}\;
\hpg21{\!\frac{\cc-k-\ell}2,\frac{1+\cc-k+\ell}{2}}{\frac12-k}{\,x}+ \\ & & 
\frac{(-1)^{k}\left(\frac{1+\cc-k-\ell}2\right)_k
\left(\frac{\cc-k-\ell}2\right)_{k+\ell+1}}
{\left(\frac{1+\cc-k-\ell}2\right)_\ell\;\left(\frac12\right)_{k+1}}\;
x^{k+\frac12}\hpg21{\!\frac{\cc-k-\ell}2,\frac{1+\cc-k+\ell}{2}}{\frac32+k}{\,x},\\
\end{eqnarray*}
We have, still with $\cc=-n/m$:
\begin{eqnarray*}
W_1\equal
\frac{(-1)^{k+\ell}\left(1-\cc\right)_{k+\ell}}{2^{k+\ell}\left(\frac{1+\cc-k-\ell}2\right)_\ell}
x^{k/2} \left(1+\sqrt{x}\right)^{-\cc}
\app3{k+1,\ell+1; -k,-\ell}{1-\cc}{\frac{\sqrt{x}+1}{2\sqrt{x}},\frac{1+\sqrt{x}}2}\\
\equal \frac{\left(\frac12\right)_k\left(\frac12\right)_\ell}{\left(\frac{1+\cc-k-\ell}2\right)_\ell}
\left(1+\sqrt{x}\right)^{-\cc+k+\ell}
\app2{\cc-k-\ell; -k,-\ell}{-2k,\,-2\ell}{\frac{2\sqrt{x}}{\sqrt{x}+1},\frac2{1+\sqrt{x}}}.
\end{eqnarray*}
The function $W_2$ is the conjugate $\sqrt{x}\mapsto-\sqrt{x}$ of $W_1$.
The polynomials $\Theta_1(x)$ and $\Theta_2(x)$, defined up to a constant multiple, 
appear in the semi-invariants
\begin{equation}
\Theta_1(x)=W_1^m+W_2^m, \qquad  x^{k+\frac12}\,\Theta_2(x)=W_1^m-W_2^m.
\end{equation}
The semi-invariants of degree 2 and $m$ are related by formula (\ref{eq:insrel}), with
$C=4^{n+1}\left/\left(\frac{1-k-\ell}2-\frac{n}{2m}\right)_\ell^{2m}\right.$. 
\end{remark}

\subsection{Explicit cases}

Here we take a look at explicit expressions of Klein's pull-back coverings
for small values of $k$, $\ell$.
Of course,  for $k=0$, $\ell=0$ we have the pull-back transformations
of general dihedral functions from Section \ref{sec:trihedr}.

With $k=1$, $\ell=0$ we have hypergeometric equations $E(1/2,3/2,n/m)$.  We have to expand in (\ref{eq:gmpower}) the following expression: 
\begin{equation} \label{eq:inv1232}
(1+\sqrt{x})^n\left(1-\frac{n\,\sqrt{x}}{m}\right)^m
\end{equation}
For $n=1$ we get this explicit expansion:
\begin{equation} \label{eq:id32q}
\left(1+\sqrt{x}\right)\left(1-\frac{\sqrt{x}}{m}\right)^m=
\theta_3\!\left(x\right)+x^{3/2}\,\theta_4\!\left(x\right),
\end{equation}
where
\begin{eqnarray*}
\theta_3(x)=\hpg21{-\frac{m}2,-\frac{m+1}2}{-1/2}{\frac{x}{m^2}},\quad
\theta_4(x)=\frac{m^2-1}{3m^2}\,\hpg21{-\frac{m-2}2,-\frac{m-3}2}{5/2}{\frac{x}{m^2}}.
\end{eqnarray*}
Formula (\ref{eq:insrel}) becomes
\begin{equation}
\theta_3(x)^2-x^3\,\theta_4(x)^2=(1-x)\left(1-\frac{x}{m^2}\right)^m.
\end{equation}
This identity is a special case of Lemma \ref{th:symsq}, with $a=-m-1$, $k=1$, $z=x/m^2$. 

Expressions for $\Theta_1(x)$, $\Theta_2(x)$ appear to be more complicated for higher $n$,
even if we keep $k=1$, $\ell=0$. For $n=2$ we still can find the explicit expressions
\begin{eqnarray*}
\theta_5(x)\equal \hpg32{\!-\frac{m+1}2,-\frac{m+2}2,-\frac{m}{m+2}}
{-\frac12,-\frac{2(m+1)}{m+2}}{\frac{4x}{m^2}},\nonumber\\
\theta_6(x)\equal \frac{2(m^2\!-\!4)}{3m^2} 
\,\hpg32{\!-\frac{m-1}2,-\frac{m-2}2,\frac{m+6}{2(m+2)}}
{\frac52,\,-\frac{m-2}{2(m+2)}}{\frac{4x}{m^2}}
\end{eqnarray*}
for
\begin{equation}
{ (1+\sqrt{x})^2\!\left(1-\frac{2\sqrt{x}}m\right)^m}=
\theta_5\!\left(x\right)+x^{3/2}\,\theta_6\!\left(x\right).
\end{equation}

Let $G_{k,\ell}(\sqrt{x})$ denote the $G(\sqrt{x})$ polynomial in (\ref{eq:gmpower}) 
for particular values of $k,\ell$. We have
\begin{eqnarray*}
G_{2,0}(\sqrt{x}) \equal 1-\frac{n\,\sqrt{x}}{m}+\frac{n^2-m^2}{3m^2}\,x,\\
G_{3,0}(\sqrt{x}) \equal 1-\frac{n\,\sqrt{x}}{m}+\frac{2n^2-3m^2}{5m^2}\,x
-\frac{n(n^2-4m^2)}{15m^3}x^{3/2},\\
G_{4,0}(\sqrt{x}) \equal  1-\frac{n\,\sqrt{x}}{m}+\frac{3(n^2-2m^2)}{7m^2}\,x
-\frac{n(2n^2-11m^2)}{21m^3}x^{3/2}+\frac{(n^2-m^2)(n^2-9m^2)}{105m^3}\,x^2,\\
G_{1,1}(\sqrt{x}) \equal  1-\frac{n\,\sqrt{x}}{m}+x,\\
G_{2,1}(\sqrt{x}) \equal  1-\frac{n\,\sqrt{x}}{m}+\frac{n^2\,x}{3m^2}-\frac{n\,x^{3/2}}{3m},\\
G_{3,1}(\sqrt{x}) \equal  1-\frac{n\,\sqrt{x}}{m}+\frac{2(n^2-m^2)}{5m^2}\,x
-\frac{n(n^2-m^2)}{15m^3}\,x^{3/2}+\frac{(n^2-m^2)}{15m^4}\,x^2,\\
G_{2,2}(\sqrt{x}) \equal  1-\frac{n\,\sqrt{x}}{m}+\frac{n^2+2m^2}{3m^2}\,x-\frac{n\,x^{3/2}}{m}+x^2.
\end{eqnarray*}
The expressions $G_{k,\ell}(\sqrt{x})^m$ 
have peculiar simultaneous ``approximation" properties at $x=0$ and $x=\infty$ as well.
Specifically, the expression \mbox{$(1+\sqrt{x})^{n}\,G_{k,\ell}(\sqrt{x})^m$} has the following 
leading and trailing coefficients  zero:
\begin{itemize} 
\item At $x=0$, the coefficients to $x^{1/2}, x^{3/2}, \ldots, x^{k-1/2}$ are equal to zero;
\item At $x=\infty$, the coefficients to $x^{(N-1)/2}, x^{(N-3)/2}, \ldots, x^{(N+1)/2-\ell}$ are equal to zero;
here $N=n+m(k+\ell)$.
\end{itemize}

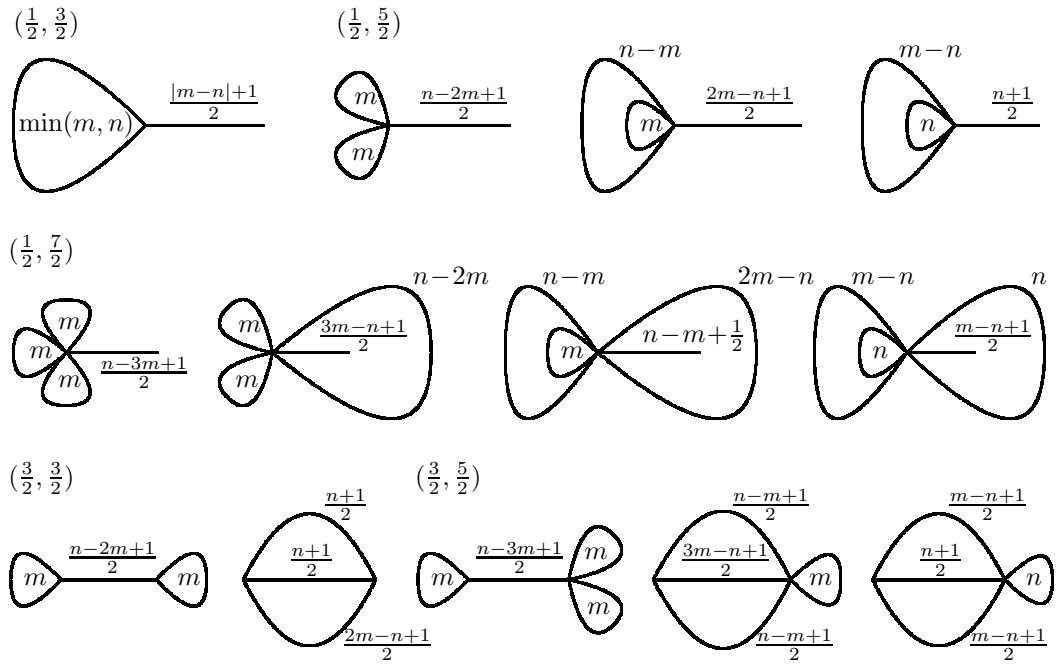
\begin{figure}
\begin{picture}(358,250)(0,-10) \thicklines
\put(-12,242){$(\frac12,\frac32)$} 
\put(-10,204){$\min(m,n)$} 
\put(46,212){$\frac{|m-n|+1}2$} \qbezier(38,206)(-12,256)(-12,206) \qbezier(38,206)(-12,156)(-12,206) \put(38,206){\line(1,0){45}} 
\put(110,242){$(\frac12,\frac52)$} 
\put(141,212){$\frac{n-2m+1}2$} 
\put(117,214){$m$} 
\put(116,193){$m$} \qbezier(130,206)(102,212)(113,223) \qbezier(130,206)(124,234)(113,223) \qbezier(130,206)(102,200)(113,189) \qbezier(130,206)(124,178)(113,189) \put(130,206){\line(1,0){46}} 
\put(225,204){$m$} 
\put(217,232){$n\!-\!m$} 
\put(249,212){$\frac{2m-n+1}2$} \qbezier(238,206)(220,224)(220,206) \qbezier(238,206)(220,188)(220,206) \qbezier(238,206)(203,256)(203,206) \qbezier(238,206)(203,156)(203,206) \put(238,206){\line(1,0){48}} 
\put(331,204){$n$} 
\put(323,232){$m\!-\!n$} 
\put(357,212){$\frac{n+1}2$} \qbezier(344,206)(326,224)(326,206) \qbezier(344,206)(326,188)(326,206) \qbezier(344,206)(309,256)(309,206) \qbezier(344,206)(309,156)(309,206) \put(344,206){\line(1,0){32}} 
\put(-14,156){$(\frac12,\frac72)$} 
\put(-6,118){$m$} 
\put(20,110){$\frac{n-3m+1}2$} 
\put(5,129){$m$} 
\put(5,107){$m$} \qbezier(8,120)(-12,138)(-12,120) \qbezier(8,120)(-12,102)(-12,120) \qbezier(8,120)(26,140)(8,140) \qbezier(8,120)(-10,140)(8,140) \qbezier(8,120)(26,100)(8,100) \qbezier(8,120)(-10,100)(8,100) \put(8,120){\line(1,0){35}} 
\put(103,124){$\frac{3m-n+1}2$} 
\put(73,128){$m$} 
\put(72,107){$m$} 
\put(139,146){$n\!-\!2m$} \qbezier(86,120)(146,170)(146,120) \qbezier(86,120)(146,70)(146,120) \qbezier(86,120)(58,126)(69,137) \qbezier(86,120)(80,148)(69,137) \qbezier(86,120)(58,114)(69,103) \qbezier(86,120)(80,92)(69,103) \put(86,120){\line(1,0){29}} 
\put(195,118){$m$} 
\put(188,146){$n\!-\!m$} 
\put(226,124){$n\!-\!m\!+\!\frac{1}2$} 
\put(262,146){$2m\!-\!n$} \qbezier(209,120)(190,138)(190,120) \qbezier(209,120)(190,102)(190,120) \qbezier(209,120)(174,170)(174,120) \qbezier(209,120)(174,70)(174,120) \qbezier(209,120)(269,170)(269,120) \qbezier(209,120)(269,70)(269,120) \put(209,120){\line(1,0){39}} 
\put(313,118){$n$} 
\put(305,146){$m\!-\!n$} 
\put(343,124){$\frac{m-n+1}2$} 
\put(373,146){$n$} \qbezier(326,120)(308,138)(308,120) \qbezier(326,120)(308,102)(308,120) \qbezier(326,120)(291,170)(291,120) \qbezier(326,120)(291,70)(291,120) \qbezier(326,120)(378,170)(378,120) \qbezier(326,120)(378,70)(378,120) \put(326,120){\line(1,0){26}} 
\put(-14,70){$(\frac32,\frac32)$} 
\put(-8,32){$m$} 
\put(8,40){$\frac{n-2m+1}2$} 
\put(50,32){$m$} \qbezier(6,34)(-13,54)(-13,34) \qbezier(6,34)(-13,14)(-13,34) \put(6,34){\line(1,0){36}} \qbezier(42,34)(61,54)(61,34) \qbezier(42,34)(61,14)(61,34) 
\put(92,39){$\frac{n+1}2$} 
\put(105,60){$\frac{n+1}2$} 
\put(112,7){$\frac{2m-n+1}2$} \put(75,34){\line(1,0){50}} \qbezier(125,34)(100,84)(75,34) \qbezier(125,34)(100,-16)(75,34) 
\put(140,70){$(\frac32,\frac52)$} 
\put(146,32){$m$} 
\put(162,40){$\frac{n-3m+1}2$} 
\put(204,42){$m$} 
\put(205,21){$m$} \qbezier(160,34)(141,54)(141,34) \qbezier(160,34)(141,14)(141,34) \put(160,34){\line(1,0){38}} \qbezier(198,34)(226,40)(215,51) \qbezier(198,34)(204,62)(215,51) \qbezier(198,34)(226,28)(215,17) \qbezier(198,34)(204,6)(215,17) 
\put(289,32){$m$} 
\put(240,39){$\frac{3m-n+1}2$} 
\put(259,60){$\frac{n-m+1}2$} 
\put(268,7){$\frac{n-m+1}2$} \put(230,34){\line(1,0){52}} \qbezier(282,34)(301,53)(301,34) \qbezier(282,34)(301,15)(301,34) \qbezier(282,34)(257,86)(230,34) \qbezier(282,34)(257,-18)(230,34) 
\put(371,32){$n$} 
\put(330,39){$\frac{n+1}2$} 
\put(341,60){$\frac{m-n+1}2$} 
\put(349,7){$\frac{m-n+1}2$} \put(313,34){\line(1,0){50}} \qbezier(363,34)(381,52)(381,34) \qbezier(363,34)(381,16)(381,34) \qbezier(363,34)(338,84)(313,34) \qbezier(363,34)(338,-16)(313,34) 
\end{picture}
\caption{{\it Dessins d'enfant} of Klein pullback coverings for dihedral hypergeometric equations, schematically}
\centering \label{fig:dessins}
\end{figure}

As typical for hypergeometric transformations, the coverings for (\ref{eq:ledpb}) and (\ref{eq:ledpc})
are {\em Belyi maps} as they ramify only above three points. They are topologically determined
by their {\em dessins d'enfant} \cite{zvonkin}, 
which are preimages of, say, $[0,1]\subset\CC$ with the 
endpoints $z=0$, $z=1$ marked black and white. 
The {\em dessin d'enfant} for (\ref{eq:ledpb}) is a line segment with
the black and white points alternating on it 
--- which is a well-known pattern for Tchebyshev polynomials.
The {\em dessin d'enfant} for (\ref{eq:ledpc}) is a circle segment with $n$ black and $n$ 
white points alternating on it.

Figure \ref{fig:dessins} schematically depicts the {\em dessins d'enfant} for 
Klein transformations of dihedral function with $k+\ell\le 3$. Each loop with coinciding endpoints
represents a circle with alternating black and white points, and the number of each kind of points
is given inside the loop or beside it. Each segment with distinct endpoints represent
a linear segment with alternating black and white points, and the number given beside it
is the arithmetic average of the number of black and white points on it. The number is integer
precisely when the endpoints have different colour. 
The first row depicts the dessins of Klein morphisms for $E(1/2,3/2,n/m)$, and
the dessins for $E(1/2,5/2,n/m)$ for the cases $n>2m$, $m<n<2m$ and $n<m$, respectively.
The second row depicts the dessins of Klein morphisms for $E(1/2,7/2,n/m)$
for the cases $n>3m$, $2m<n<3m$, $m<n<2m$ and $n<m$, respectively.
The third row depicts the dessins for $E(3/2,3/2,n/m)$ when $n>2m$ or $n<2m$,
and the dessins for $E(3/2,5/2,n/m)$ when $n>3m$, $m<n<3m$ or $n<m$. 
As we see, the variation of the dessins is quite involved. The cases when $n$ is an integer
multiple of $m$ are either non-existent or should have a bordering type, 
depending on whether the respective hypergeometric equation has logarithmic solutions
or a trivial monodromy group, as distinguished by \cite[Theorem 5.1]{tdihedral}.

\begin{remark} \rm
Attractive trigonometric formulas in \cite[Section 6]{tdihedral} suggest that it is worthwhile
to consider trigonometric expressions for the (semi)-invariants 
$\Theta_1(x)$, $\Theta_2(x)$, $\Psi(x)$. 
In the case of $E(1/2,3/2,n/m)$, expression (\ref{eq:inv1232}) 
is equated to 
\mbox{$\Theta_1(x)+x^3\Theta_2(x)$}.
After the substitution $x\mapsto -\tan^2 t$ and multiplication by $\cos^{n+m}t$, 
expression (\ref{eq:inv1232}) becomes 
\begin{equation}
(\cos nt+i\sin nt)\left(\cos t-\frac{in}{m}\sin t\right)^{\!m}.
\end{equation}
The exponential form of this expression is
\begin{equation}
\left(\cos^2 t+\frac{n^2}{m^2}\sin^2 t\right)^{m/2}
\exp i\!\left(nt-m\arctan\left(\frac{n\tan t}{m}\right)\right).
\end{equation}
The real and complex parts of this formula basically give trigonometric forms
of the semi-invariants $\Theta_1$, $\Theta_2$. 
Up to a power of $\cos t$, both sides of (\ref{eq:insrel}) express
the squared norm $\left(\cos^2 t+\frac{n^2}{m^2}\sin^2 t\right)^{m}$.
\end{remark}

\subsection{Differential expressions}
\label{diffinv}

The components of the pull-back map $\Psi,\Theta_1,\Theta_2$,
which are also monodromy invariants, satisfy differential equations
that can be derived from necessary properties of the pull-back maps.
The results are summarised in the following theorem.

The constant $\Psi(1)$ in (\ref{eq:difftpsii}) can be found using (\ref{eq:psii1}).
Equation (\ref{eq:diffpsi}) is just a symmetric tensor square equation
for a dihedral hypergeometric equation, normalized  to have a polynomial solution.
Equation (\ref{eq:thetaa}) can be seen as a Fuchsian equation for $\Theta_2$ 
with the singularities $x=0$, $x=1$, $x=\infty$ and the roots of $\Psi$. 
The polynomial $\Psi$ can be eliminated using (\ref{eq:diffpsi})
and a differential Gr\"obner basis algorithm, but 
that gives a complicated non-linear differential equation (of order 4, algebraic degree 4)
for $\Theta_2$ alone.
\begin{theorem}
The following differential relations for $\Psi(x)$ and $\Theta_2(x)$ 
hold:
\begin{align} \label{eq:difftpsii}
x(x-1) \left(\Psi'{}^2-2\Psi''\Psi\right)+\big((2k-3)x-2k+1\big)\Psi'\Psi\hspace{120pt}\nonumber\\
=\left((k+\ell)(k-\ell-1)-\frac{n^2}{m^2(x-1)}\right) \Psi^2+\frac{n^2\,\Psi(1)^2\,x^{2k}}{m^2(x-1)},
\end{align}
\begin{align} \label{eq:diffpsi}
\frac{d^3\Psi}{dx^3}=& \,
3\left( \frac{k-\frac12}{x} +\frac{1}{x-1} \right) \frac{d^2\Psi}{dx^2} \nonumber\\
&-\left(\frac{k(2k-1)}{x^2}+\frac{m^2-n^2}{m^2(x-1)^2}
+\frac{k^2-\ell^2-5k-\ell+2+n^2/m^2}{x(x-1)} \right)\frac{d\Psi}{dx}\nonumber\\
&+\frac{(k+\ell)(k-\ell-1)(2k-1)x-2k\big((k+\ell)(k-\ell-1)+n^2/m^2\big)}{x^2(x-1)^2}\,\Psi.
\end{align}
\begin{align} \label{eq:thetaa}
\frac{\Theta_2''}{\Theta_2}
+\left(\frac{k+\frac32}{x}-\frac{n-1}{x-1}-\frac{(m-1)\Psi'}{\Psi} \right) \frac{\Theta_2'}{\Theta_2}
+\frac{m(m-1)}2\frac{\Psi''}{\Psi} \nonumber\\
-\frac{m-1}{4}\left( \frac{m(2k-1)+4k+2}{x}+\frac{2(m+n)}{x-1}\right)\frac{\Psi'}{\Psi} \nonumber\\
+\frac{m^2(k+\ell)(k-\ell-1)+(n-2)(n-2k-1)}{4x(x-1)} &=0.\qquad\qquad
\end{align}
\end{theorem}
\begin{proof} Let $d=(k+\ell)m+n$ denote the degree of the pullback covering.
By Corollary \ref{th:invsrel}, the covering has the form
\begin{equation}
\varphi(x)=\frac{x^{2k+1}\,\Theta_2(x)^2}{\Theta_1(x)^2}, \qquad
\varphi(x)-1=\frac{C\,(1-x)^n\,\Psi(x)^m}{\Theta_1(x)^2}.
\end{equation}
The degree of $\Psi(x)$ is $k+\ell$. If $d$ is odd then the point $x=\infty$ lies above $z=\infty$.
If $d$ is even, then the point $x=\infty$ lies above $z=0$ just as $x=0$.

Let $\alpha,\beta$ denote the degree 
of, respectively, $\Theta_1(x),\Theta_2(x)$. Then 
\[ \begin{array}{lll}
\alpha=\frac{d-1}2-\ell, & \beta=\frac{d-1}2-k, & \mbox{ if $d$ is odd},\\
\alpha=\frac{d}2, & \beta=\frac{d}2-k-\ell-1, & \mbox{ if $d$ is even}.
\end{array} \]
Computations are more instructive and easier to follow when we consider
an alternative normalization of the covering.
Namely, we make $x=0$ and $x=1$  the points with half-integer local exponent differences.
Hence we interchange the points $x=1$ and $x=\infty$.
We define the polynomials $P(x)$, $Q(x)$, $R(x)$ to be the monic polynomials
proportional to, respectively,
$$\textstyle
 (x-1)^{\beta}\,\Theta_2\big(\frac{x}{x-1}\big), \quad
 (x-1)^{\alpha}\,\Theta_1\big(\frac{x}{x-1}\big), \quad 
 (x-1)^{k+\ell}\,\Psi\big(\frac{x}{x-1}\big).
$$


When the degree $d$ of the pull-back covering is odd,
the pull-back covering has the form
\begin{equation}
\varphi(x)=\frac{c\,x^{2k+1}\,P^2}{R^m}, \qquad
\varphi(x)-1=\frac{c\,(x-1)^{2\ell+1}\,Q^2}{R^m}.
\end{equation}
Since we know all the ramification points, we know the numerator
of the derivative:
\begin{equation}
\varphi'(x)=\frac{ncx^{2k}(x-1)^{2\ell}PQ}{R^{m+1}}.
\end{equation}
The factor $nc$ follows from local consideration at $x=\infty$.
We have two ways to express the logarithmic derivatives of $\varphi(x)$ and $\varphi(x)-1$,
and that gives the equations
\begin{align} \label{eq:logdiffs}
\frac{n\,(x-1)^{2\ell}Q}{x\,PR} = \frac{2k+1}{x}+\frac{2P'}{P}-\frac{mR'}{R},\qquad
\frac{n\,x^{2k}\,P}{(x-1)QR} = \frac{2\ell+1}{x-1}+\frac{2Q'}{Q}-\frac{mR'}{R}.
\end{align}
Elimination of $Q$ gives
\begin{align} \label{eq:logdiffp}
n^2\,x^{2k}(x-1)^{2\ell}P =& \, 
 x(x-1)\left(4P''R^2-4(m-1)P'R'R-2mPR''R+m^2PR'{}^2\right) \nonumber \\
& +\left( 2(k-\ell+2)x-2k-3\right) \left(2P'R-mPR'\right) R \nonumber \\
& +(2k+1) \left( 2(x-1)R'-(2\ell-1) R\right) P R.
\end{align}
Additionally, the pull-backed equation can be computed symbolically.
A tedious computation shows that the pullback
$z\mapsto \varphi(x)$, $y(z)\mapsto Y(\varphi(x))\big/\sqrt{R}$
of the hypergeometric equation for $\hpg21{\!1/{2m},-1/{2m}}{1/2}{\,z}$ is
\begin{align} \label{eq:pbacked}
\frac{d^2Y}{dx^2}=&\left( \frac{k-\frac12}{x} +\frac{\ell-\frac12}{x-1} \right) \frac{dY}{dx} \nonumber\\
&+\left( \frac{R''}{2R}-\frac{R'{}^2}{4R^2}
-\left( \frac{k-\frac12}{x} +\frac{\ell-\frac12}{x-1} \right) \frac{R'}{2R}
+\frac{n^2x^{2k-1}(x-1)^{2\ell-1}}{4m^2R^2} \right) Y.
\end{align}
This equation must coincide with the hypergeometric equation for
$$\hpg21{\frac{n}{2m}-\frac{k+\ell}2,-\frac{n}{2m}-\frac{k+\ell}2}{\frac12-k}{\,x},$$
which has the form
\begin{align*}
\frac{d^2Y}{dx^2}=&\left( \frac{k-\frac12}{x} +\frac{\ell-\frac12}{x-1} \right) \frac{dY}{dx}
+\frac1{4x(x-1)} \left(\frac{n^2}{m^2}-(k+\ell)^2 \right) Y.
\end{align*}
The coefficients to $dY/dx$ coincide, but equating the coefficients to $Y$
gives a new non-linear differential equation for $R$:
\begin{align} \label{eq:difftrr}
\left((k+\ell)^2-\frac{n^2}{m^2}\right) R^2+\frac{n^2}{m^2}\,x^{2k}(x-1)^{2\ell}
=& \left(2(k+\ell-1)x-2k+1\right)R'R \nonumber\\
& +x(x-1) \left(R'{}^2-2R''R\right).
\end{align}
Application of the differential operator $d/dx-2k/x-2\ell/(x-1)$ removes the terms with 
$x^{2k}(x-1)^{2\ell}$ and $R'{}^2$, giving the Fuchsian equation
\begin{align} \label{eq:rsqeq}
\frac{d^3R}{dx^3}=& \,
3\left( \frac{k-\frac12}{x} +\frac{\ell-\frac12}{x-1} \right) \frac{d^2R}{dx^2} \nonumber\\
&-\left(\frac{k(2k-1)}{x^2}+\frac{\ell(2\ell-1)}{(x-1)^2}
+\frac{(k+\ell)^2+(2k-1)(2\ell-1)-n^2/m^2}{x(x-1)} \right)\frac{dR}{dx}\nonumber\\
&+\left((k+\ell)^2-\frac{n^2}{m^2}\right)\frac{(k+\ell)x-k}{x^2(x-1)^2}\,R.
\end{align}
This is the symmetric tensor square equation for $R$, as expected. 
Elimination of the $x^{2k}(x-1)^{2\ell}$ term from (\ref{eq:logdiffp}) and (\ref{eq:difftrr})
gives the equation
\begin{align} \label{eq:ppdiffeq}
\frac{P''}{P}+\left(\frac{k+\frac32}{x}-\frac{\ell-\frac12}{x-1}-\frac{(m-1)R'}{R} \right) \frac{P'}{P}
+\frac{m(m-1)}2\frac{R''}{R} \nonumber\\
-\frac{m-1}{4}\left( \frac{m(2k-1)+4k+2}{x}+\frac{m(2\ell-1)}{x-1}\right)\frac{R'}{R} \nonumber\\
+\frac{m^2(k+\ell)^2-n^2-(2k+1)(2\ell-1)}{4x(x-1)} &=0.
\end{align}

When the degree of the pull-back covering is even,
the pull-back covering has the form
\begin{equation}
\varphi(x)=\frac{c\,x^{2k+1}(x-1)^{2\ell+1}P^2}{R^m}, \qquad
\varphi(x)-1=\frac{c\,Q^2}{R^m},
\end{equation}
the transformed equation (\ref{eq:pbacked}) and the equations (\ref{eq:difftrr})--(\ref{eq:rsqeq})
for $R$ do not change, but the logarithmic derivative expressions (\ref{eq:logdiffs}) change to
\begin{align*}
\frac{n\,Q}{x(x-1)PR} = \frac{2k+1}{x}+\frac{2\ell+1}{x-1}+\frac{2P'}{P}-\frac{mR'}{R},\qquad
\frac{n\,x^{2k}(x-1)^{2\ell}P}{QR} = \frac{2Q'}{Q}-\frac{mR'}{R}.
\end{align*}
Consequently, equation (\ref{eq:logdiffp}) changes to
\begin{align*} 
n^2\,x^{2k}(x-1)^{2\ell}P =& \, 
 x(x-1)\left(4P''R^2-4(m-1)P'R'R-2mPR''R+m^2PR'{}^2\right) \nonumber \\
& +\left( 2(k+\ell+3)x-2k-3\right) \left(2P'R-mPR'\right) R \nonumber \\
& -2(2k+1)PR'R+4 (k+\ell+1) (xR'+R) P R.
\end{align*}
and equation (\ref{eq:ppdiffeq}) changes to
\begin{align} \label{eq:ppdiffeq2}
\frac{P''}{P}+\left(\frac{k+\frac32}{x}+\frac{\ell+\frac32}{x-1}-\frac{(m-1)R'}{R} \right) \frac{P'}{P}
+\frac{m(m-1)}2\frac{R''}{R} \nonumber\\
-\frac{m-1}{4}\left( \frac{m(2k-1)+4k+2}{x}+\frac{m(2\ell-1)+4\ell+2}{x-1}\right)\frac{R'}{R} \nonumber\\
+\frac{m^2(k+\ell)^2-n^2+4(k+\ell+1)}{4x(x-1)} &=0.
\end{align}

The transformation $x\mapsto x/(x-1)$ back to $\Psi(x),\Theta_2(x)$
transforms (\ref{eq:difftrr}) to (\ref{eq:difftpsii}), equation (\ref{eq:rsqeq}) to (\ref{eq:diffpsi}),
and both equations (\ref{eq:ppdiffeq}) and (\ref{eq:ppdiffeq2}) to (\ref{eq:thetaa}). 
\end{proof}

\subsection{All pullback transformations}

Theorem \ref{th:genftrans} states that Klein pull-back transformations
of dihedral hypergeometric equations are unique up to M\"obius transformations.
Other than that, there are:
\begin{itemize}
\item Specializations of classical transformations \cite{goursat} for hypergeometric functions.
It is straightforward to check the entries of \cite[Table 1]{algtgauss}
whether their specializations involve dihedral functions.
Here is the list of relevant classical transformations (except the quadratic):
\begin{eqnarray*} 
\textstyle \hspace{-20pt}E\left(\frac12,\frac13,\frac{2k+1}4\right) \stackrel{3}{\longleftarrow} 
E\left(\frac12,k+\frac12,\frac{2k+1}4\right),& \quad &
\textstyle E\left(\frac12,\frac13,k+\frac12\right)\stackrel{3}{\longleftarrow}
E\left(\frac12,k+\frac12,2k+1\right),\\
\textstyle E\left(\frac12,\frac12,\frac13\right)\stackrel{4}{\longleftarrow}
E\left(\frac12,\frac32,\frac13\right),&&
\textstyle E\left(\frac12,\frac12,\frac14\right)\stackrel{4}{\longleftarrow}
E\left(\frac12,\frac12,1\right),\\
\textstyle E\left(\frac12,\frac13,\frac13\right)\stackrel{3}{\longleftarrow}
E\left(\frac12,\frac12,\frac12\right),&&
\textstyle E\left(\frac12,\frac12,\frac13\right)\stackrel{6}{\longleftarrow}
E\left(\frac12,\frac12,2\right),\\
\textstyle E\left(\frac12,\frac13,\frac14\right)\stackrel{6}{\longleftarrow}
E\left(\frac12,\frac12,\frac12\right),&&
\textstyle E\left(\frac12,\frac12,\frac13\right)\stackrel{6}{\longleftarrow} E(1,1,1).
\end{eqnarray*}
\item Transformations (\ref{eq:ledpb}), which are unique as well;
\item Transformations from $E(1/2,1/2,1/m)$ to hypergeometric equations
with a smaller dihedral group or a cyclic monodromy group. 
\end{itemize}
Transformations to smaller dihedral groups are unique, 
since the Schwartz maps are identified similarly uniquely.
Even the formulas of Section \ref{adihedral} are correct when $\gcd(n,m)>1$
but $n/m\not\in\ZZ$.  

But uniqueness is not guaranteed for transformations from $E(1/2,1/2,1/m)$
to a cyclic monodromy group, even if the branching pattern is fixed.
For a counterexample, consider coverings with the branching pattern
$2+2+2+2+2=4+2+2+1+1=5+5$, suitable for the transformation
$E(1/2,1/2,1/5)\stackrel{\,10}{\longleftarrow}E(1/2,1/2,2)$.
An explicit covering is
\begin{equation}
f_{10}(x)= \frac{25(5-\sqrt{5})\,x\,(x-1)\left(x^2-x+\frac{2-\sqrt5}{16}\right)^2}
{2048 \left(x^2-x+\frac{3-\sqrt5}{32}\right)^5}.
\end{equation}
The conjugation $\sqrt5\mapsto-\sqrt5$ gives a different covering, since the only
candidate M\"obius transformation to realize the symmetry $x\mapsto1-x$
leaves the covering intact. Correspondingly, there are also two {\em dessins d'enfant}
with the same branching data:
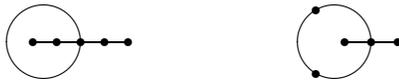
\[
\begin{picture}(190,32)(0,-16)
\put(40,0){\circle{26}}    \put(150,0){\circle{26}}
\put(36,0){\line(1,0){35}}   \put(154,0){\line(1,0){20}}
\put(36,0){\circle*{3}} \put(45,0){\circle*{3}}  \put(54,0){\circle*{3}}  
\put(63,0){\circle*{3}} \put(72,0){\circle*{3}}
\put(154,0){\circle*{3}} \put(164,0){\circle*{3}}  \put(174,0){\circle*{3}}  
\put(143,12){\circle*{3}} \put(143,-12){\circle*{3}}
\end{picture}
\]
In fact, there is a continuous family of transformations 
$E(1/2,1/2,1/5)\stackrel{\,10}{\longleftarrow}E(1/2,1/2,2)$
with a varying location for the $x$-point with the local exponent difference 2.
This family is a composition of the transformation 
$E(1/2,1/2,1/5)\stackrel{\,5}{\longleftarrow} E(1/2,1/2,1)$
with a continuous family of quadratic transformations  
$E(1/2,1/2,1)\stackrel{\,2}{\longleftarrow} E(1/2,1/2,2)$. 
The degree 5 covering is realized by the rational function
\[
f_5(z)=\frac{z(z^2+10z+5)^2}{(z-1)^5},
\]
while the variable quadratic covering can be taken to be
$z(x)=c\,x(1-x)/\left(x-\frac12\right)^2$ so that its branching points
are $x=1/2$ and $x=\infty$. The point with the local exponent difference is always $x=\infty$,
and it's projection $z=c$ is a branching point in the fiber $f_5(z)=0$ when
$c=5\pm2\sqrt5$ giving two composite coverings with the required branching pattern.
For $c=(5\pm2\sqrt5)/5$ we get two different coverings with the branching pattern
$4+2+2+2=2+2+2+2+1+1=5+5$.
The non-uniqeness comes from the freedom of choosing a Schwarz map
for $E(1/2,1/2,2)$. Following the identification of Schwarz maps
around formulas (\ref{kleing0}) and (\ref{eq:kleing2}), we have
a unique (up to a constant multiple) solution $x^2\,\hpg21{\!1,\,3/2}{3}{x}$
with the local exponent 2, but a continuous family of ''hypergeometric" solutions
$\hpg21{\!-1,-1/2}{-1}{x}$ with the local exponent 0. 

The claim in \cite[Remark 7.1]{algtgauss} that transformations from Klein's standard 
hypergeometric equations with algebraic solutions are unique (up to M\"obius transformations)
thus appears to be wrong when the transformed equation has a cyclic monodromy group.
This is also a case when non-Belyi maps appear as pull-back coverings.

\begin{remark} \rm \label{rm:degtrans}
As observeded in \cite[Remark 4.1, Subsection 5.4]{tdihedral}, hypergeometric equations
with terminating or logarithmic solutions have a degenerate structure of the 24 Kummer's
solutions of a general hypergeometric equation. If a pull-back transformation relates 
hypergeometric equations with a degenerate structure (at least on one side),
then standard relations between solutions of the hypergeometric equations may not hold,
as exemplified in \cite[Remark 4.1]{tdihedral}. In transformation (\ref{eq:ledpc}), 
the set of 24 Kummer's solutions on the $(1,n\cc,n\cc)$ side is always degenerate.
A special care must be taken to relate hypergeometric functions under the transformations
of Sections \ref{sec:trihedr} and \ref{adihedral} when terminating or logarithmic solutions are involved.
But the transformations always hold as pullbacks of hypergeometric equations.

Transformations  (\ref{eq:quatriv}) and (\ref{eq:quadih}) involve degenerate Gauss hypergeometric functions on both sides when $\cc$ is an integer. Let us assume here that $\cc$ is an integer. 
By \cite[Theorem 5.1]{tdihedral}, 
logarithmic solutions appear in (\ref{eq:quadih}) or  (\ref{eq:quatriv})
if and only if $|\cc|\le k$. The special case $k=0$ of (\ref{eq:quadih}) leads to Tchebyshev polynomials, 
as demonstrated in (\ref{eq:tchebysh}).
More generally,  Gegenbauer polynomials $C^{\gamma}_n$ with $\gamma=k+1$ appear 
in (\ref{eq:quadih}) when $\cc=n+k+1$. 
The special case $\cc=0$ of (\ref{eq:quatriv}) leads to Legendre polynomials.
More generally, Gegenbauer polynomials $C^{\gamma}_n$ with half-integer 
$\gamma=\frac12+|\cc|$ appear when $k=n+|\cc|$. 
For a complete picture, 
a quadratic transformation for general Gegenbauer polynomials transforms 
their hypergeometric equations as follows:
\begin{eqnarray}
& E\left(\frac12,\,\frac12-\gamma,\,n+\gamma\right)\stackrel{2}{\longleftarrow}
E\left( \frac12-\gamma,\,\frac12-\gamma,\,2n+2\gamma\right).
\end{eqnarray}
\end{remark}

\begin{remark} \rm
An other question of \cite[Remark 7.1]{algtgauss} is existence of pull-back transformations
of hypergeometric equations that do not yield a two-term identity between their solutions.
This question can be answered positively by considering a composition of
$E(1/2,1/2,1/2)\stackrel{\,4}{\longleftarrow}E(1,1,1)$ and a general
transformation $E(1,1,1)\stackrel{\,3}{\longleftarrow}E(3,2,2)$, or any other general
Klein transformation for a trivial monodromy group. 
Particularly, we take the degree 4 covering to be $-4\zeta^2/(\zeta^2-1)^2$, 
where $\zeta(z)$ is a generic M\"obius transformation, and the cubic covering
to be $z(x)=x^2 (3-2x)$. Then hypergeometric solutions of $E(1/2,1/2,1/2)$ pull-back
to general linear functions as solutions of the $E(1,1,1)$ equation, and cannot be identified
with a few ``hypergeometric" solutions of $E(1,1,1)$ or $E(3,2,2)$.
\end{remark}

\bibliographystyle{plain}
\bibliography{../../hypergeometric}

\end{document}